\documentclass[a4paper,reqno,10pt]{amsart}

\usepackage{amsthm,amsfonts,amsxtra,amssymb,amscd}

\usepackage[all]{xy}

\usepackage{enumerate}

\theoremstyle{definition}
\newtheorem{definition}{Definition}[section]

\newtheorem{example1}[definition]{Example}

\newtheorem{remark1}[definition]{Remark}

\theoremstyle{plain}

\newtheorem{corollary}[definition]{Corollary}

\newtheorem{proposition}[definition]{Proposition}
\newtheorem{lemma}[definition]{Lemma}

\newtheorem{cor}[definition]{Corollary}
\newtheorem{prp}[definition]{Proposition}
\newtheorem{lem}[definition]{Lemma}
\newtheorem*{theoremA}{Theorem A}
\newtheorem*{theoremB}{Theorem B}

\newcommand{\mN}{\mathbb N}  \newcommand{\mP}{\mathbb P} 
 \newcommand{\mR}{\mathbb R}

\newcommand{\mZ}{\mathbb Z}

 \newcommand{\calF}{\mathcal F} 
 \newcommand{\calI}{\mathcal I} \newcommand{\calJ}{\mathcal J}
 \newcommand{\calL}{\mathcal L} 
 \newcommand{\calO}{\mathcal O} 
  \newcommand{\calS}{\mathcal S}

  \newcommand{\gog}{\mathfrak g}
\newcommand{\goh}{\mathfrak h}  
 \newcommand{\gol}{\mathfrak l} \newcommand{\gom}{\mathfrak m}

\newcommand{\got}{\mathfrak t} \newcommand{\gou}{\mathfrak u}

\newcommand{\sfA}  {\mathsf A} 
\newcommand{\sfB}{\mathsf B} \newcommand{\sfC}{\mathsf C} \newcommand{\sfD}{\mathsf D} 
\newcommand{\sfE}{\mathsf E} \newcommand{\sfF}{\mathsf F} \newcommand{\sfG}{\mathsf G}

 \newcommand{\mbU}{\mathbf U}

\newcommand{\gra}{\alpha} \newcommand{\grb}{\beta}   \newcommand{\grg}{\gamma}
 
 \newcommand{\grl}  {\lambda}
\newcommand{\grs}  {\sigma}
\newcommand{\grf}  {\varphi}
  \newcommand{\grD}  {\Delta}
\newcommand{\grO}  {\Omega} \newcommand{\grL}  {\Lambda}

\newcommand{\mi}  {\imath}
\newcommand{\mj}  {\jmath}
\newcommand{\mk}  {\Bbbk}

\newcommand{\bra}        {\langle}
\newcommand{\ket}        {\rangle}
\newcommand{\End}        {\text{End}}

\renewcommand{\Im}       {\operatorname{Im}}

\newcommand{\id}      {\mathrm{id}} 

\newcommand{\incluso} {\hookrightarrow}

\newcommand{\lra}     {\longrightarrow}
\newcommand{\isocan}  {\simeq}
\newcommand{\vuoto}   {\varnothing}

\newcommand{\cech}    {\spcheck}

            \newcommand{\st}       {\, | \,}
       
         \newcommand{\mand}     {\text{ and }}
        \newcommand{\mif}      {\text{ if }}
  
\newcommand{\then}    {\Rightarrow}        

\newcommand{\wt}[1]      { {\widetilde {#1} } }
\newcommand{\wbar}[1]    { {\overline  {#1} } }

\newcommand{\iacc}   {\`{\i} }

\newcommand{\N}{\mathbb{N}}
\newcommand{\Z}{\mathbb{Z}}

\newcommand{\R}{\mathbb{R}}

\newcommand{\al}{\alpha}
\newcommand{\be}{\beta}
\newcommand{\ga}{\gamma}
\newcommand{\de}{\delta}
\newcommand{\si}{\sigma}

\newcommand{\la}{\lambda}

\newcommand{\om}{\omega}
\newcommand{\De}{\Delta}
\newcommand{\La}{\Lambda}
\newcommand{\Ga}{\Gamma}
\newcommand{\Om}{\Omega}

\newcommand{\coal}{\gra^\vee}

\newcommand{\tal}{{\wt{\gra}}}

\newcommand{\osi}{\overline{\grs}}

\newcommand{\tom}{\wt{\om}}
\newcommand{\tPhi}{{\wt{\Phi}}}
\newcommand{\tDe}{\wt{\De}}

\newcommand{\tlb}{\mathcal{O}(1)}
\newcommand{\lb}{\mathcal{L}}

\newcommand{\Pj}{\mathbb{P}}
\newcommand{\field}{\mathop{\Bbbk}}
\newcommand{\lig}{\mathfrak{g}}

\newcommand{\lit}{\mathfrak{t}}
\newcommand{\supp}{\mathop{\rm supp}\nolimits}

\newcommand{\Pic}{\mathop{\rm Pic}\nolimits}

\newcommand{\height}{\mathop{\rm ht}\nolimits}
\newcommand{\hts}{\height_0^+}

\newcommand{\typeA}{\sfA}
\newcommand{\typeB}{\sfB}
\newcommand{\typeBC}{\mathsf{BC}}
\newcommand{\typeC}{\sfC}
\newcommand{\typeD}{\sfD}
\newcommand{\typeE}{\sfE}
\newcommand{\typeF}{\sfF}
\newcommand{\typeG}{\sfG}

\newcommand{\typeDL}{\typeD_5^L}
\newcommand{\typeDR}{\typeD_5^R}
\newcommand{\tpA}{\typeA}
\newcommand{\tpB}{\typeB}

\newcommand{\tpC}{\typeC}
\newcommand{\tpD}{\typeD}
\newcommand{\tpE}{\typeE}
\newcommand{\tpF}{\typeF}
\newcommand{\tpG}{\typeG}
\newcommand{\tpDL}{\typeDL}
\newcommand{\tpDR}{\typeDR}

\newcommand{\oG}{\overline{G}}
\newcommand{\oP}{\overline{P}}
\newcommand{\oH}{\overline{H}}

\newcommand{\low}{low\ }
\newcommand{\Low}{Low\ }

\newcommand{\wtI}{\calI}  
\newcommand{\ntI}{I}
\newcommand{\wtJ}{\calJ}  
\newcommand{\ntJ}{J}
\newcommand{\wtS}{\calS}

\newcommand{\suppw}{\supp_{\Omega}}

\newcommand{\stab}{\operatorname{stab}}
\newcommand{\senza}{\smallsetminus}
\newcommand{\retrad}{{\wt{R}}}

\newcommand{\retradJ}{{\wt{R}_\wtJ}}

\newcommand{\cvdqui}{}

\begin{document}

\title{Projective normality of complete symmetric varieties}



\author{Rocco Chiriv\iacc}

\author{Andrea Maffei}

\begin{abstract}
We prove that in characteristic zero the multiplication of sections of 
dominant line bundles on a complete symmetric variety $X=\overline{G/H}$ is a 
surjective map. As a consequence the cone defined by a complete linear system 
over $X$, or over a closed $G$ stable subvariety of $X$ is normal. 
This gives an affirmative answer to a question raised by Faltings 
in \cite{F}. A crucial point of the proof is a combinatorial property of root systems.
\end{abstract}

\maketitle


\section*{Introduction}

Let $\oG$ be an adjoint semisimple group over an algebraically closed field of characteristic zero. Given 
an involutorial automorphism $\si:\oG\rightarrow\oG$, denote by $\oH$ the subgroup of fixed points of $\si$. A wonderful compactification $X$ of the symmetric variety $\oG/\oH$ has been constructed by De Concini and Procesi \cite{CP} in characteristic zero and by De Concini and Springer \cite{CS} in positive characteristic. They describe the Picard group of $X$ as a subgroup of $\Pic(\oG/\oP)$, where $\oP$ is a suitable parabolic subgroup related to the action of $\si$ and $\oG/\oP$ is the unique closed orbit of $X$. In particular we say that a line bundle $\lb$ is dominant if $\Ga(\oG/\oP,\lb|_{\oG/\oP})\neq0$. The main result of our paper can be stated as
\begin{theoremA}
If $\lb$ and $\lb'$ are dominant line bundles on $X$, then the multiplication
$\Gamma (X,\lb)\otimes \Gamma (X,\lb')\rightarrow \Gamma (X,\lb\otimes\lb')$ is surjective.
\end{theoremA}
The projective normality of $X$ follows with a standard argument. 
Hence we give an affirmative answer to a problem raised by Faltings in \cite{F}. 
Our result has already been proved in \cite{K} in the special case of the compactification of a group, i.e. $\si:\oG\times\oG\rightarrow\oG\times\oG$, $\si(g_1,g_2)=(g_2,g_1)$, by a complete different method which does not generalize to this situation. We stress that 
it is necessary to assume that the line bundles $\lb$ and $\lb'$ are dominant as the example after the proof of Theorem~A in Section \ref{sprojectivenormality} shows.
Moreover, as De Concini explained to us, in positive characteristic the
multiplication map need not to be surjective as follows by Bruns
\cite{BR} section 4 (see also \cite{C}).

Now we briefly describe the lines of the proof of Theorem~A. We divide the $G$ modules appearing in $\Gamma (X,\lb\otimes\lb')$ in three classes and we use different strategies for each class to show that the $G$ modules appear in the image of the multiplication.
The first class is that of the modules appearing in a product of line bundles that, with respect to the dominant order, are less than $\lb$ or $\lb'$.
These are easily covered by induction on the dominant order on line bundles. The second class is formed by the modules that do not vanish when restricted to some $G$ stable subvariety $Y$ of $X$.
Notice that $Y$ is a fibration over a partial flag variety with fiber isomorphic to a complete symmetric variety $F$ with $\dim F<\dim X$ and
we can suppose that the multiplication map is surjective for $F$, using induction on dimension. In Proposition \ref{prp:ind} we show that from the surjectivity of the multiplication map on the fiber we can deduce that the multiplication for $Y$ is surjective.
So also this kind of modules appear in the image. We parametrize the remaining modules, forming the third class, introducing the notion of \low triple. Thank to the result of Section \ref{slowtriple} in which we study such triples, we can prove the surjectivity of the multiplication map for this class by a direct argument.

Our classification of \low triples is purely combinatorial and make sense for any root system $\Phi$. Suppose we have chosen a base $\De$ and let $\La^+$ be the corresponding cone of dominant weights for $\Phi$. If $(\la,\mu,\nu)$ is a triple of such weights we say that it is a \emph{\low triple} if the following conditions hold. First we require that if $\la',\mu'$ are dominant weights such that $\la'\leq\la$, $\mu'\leq\mu$ and $\nu\leq\la'+\mu'$ then $\la'=\la$ and $\mu'=\mu$. The second condition is that $\nu+\sum_{\al\in\De}\al\leq\la+\mu$. Then, if $w_0$ is the longest element of the Weyl group of $\Phi$, the result in Section \ref{slowtriple} is
\begin{theoremB}
The triple $(\la,\mu,\nu)$ of dominant weights is a \low triple if and only if $\la$ and $\mu$ are minuscule weights, $\mu=-w_0\la$ and $\nu=0$.
\end{theoremB}
Our proof of this theorem is somewhat unsatisfactory. Although we succeeded in developing a bit of general treatment, reducing in a significant way the computations involved, in some steps we still have to use a case by case analysis.

\vskip .5cm We would like to thank Corrado De Concini for useful 
conversations. We would like also to thank Paolo Papi for bringing 
to our attention the work of Stembridge \cite{Stem} and Jan Draisma 
and Arjeh Cohen for helping us coding a LiE program to check Theorem~B 
for exceptional Weyl groups in an early stage of development of this paper.


\section{Recalls on complete symmetric varieties}\label{srecalls}

In this section we collect some preliminary results for the sequel setting up
notation and reviewing the construction of the wonderful
compactification of $G/H$ (for details see \cite{CP} and \cite{CS}).

Let $\gog$ be an adjoint semisimple Lie algebra over an algebraically closed 
field $\field$ of
characteristic zero, and let $\si$ be an involutorial automorphism of $\gog$.
Denote by $\goh$ the subalgebra of fixed points of $\si$ in $\gog$. 
If $\got$ is a $\si$ stable
toral subalgebra of $\gog$, we can decompose $\lit$ as
$\lit_0\oplus\lit_1$ with $\lit_0$ the $(+1)$ eigenspace of $\si$ and $\lit_1$
the $(-1)$ eigenspace. Recall that any $\si$ stable
toral subalgebra of $\gog$ is contained in a maximal one which is itself $\si$ stable.
We fix such a $\si$ stable maximal toral subalgebra $\got$ 
for which $\dim \got_1$ is maximal
and denote this dimension by $\ell$.

Let $\Phi\subset\lit^*$ be the root system of $\lig$ 
and let $\gog = \got \oplus \bigoplus_{\gra \in \Phi} \gog_\gra$ be the root space decomposition with respect to the action of $\got$. Observe that
$\grs$ acts also on $\got^*$ and that it preserves $\Phi$ and the 
Killing
form $(\cdot,\cdot)$ on $\lit$ and $\lit^*$.
Let $\Phi_0=\{\al\in\Phi\ |\ \si(\al)=\al\}$
and $\Phi_1=\Phi\senza\Phi_0$. 
The choice of a $\grs$ stable toral subalgebra for which $\dim \got_1$ is
maximal is equivalent to the condition $\grs|_{\gog_\gra} = \id|_{\gog_\gra}$ for
all $\gra \in \Phi_0$. Moreover we can choose the set of positive roots $\Phi^+$ in such a way that $\si(\al)\in\Phi^-$ for all roots $\al\in\Phi^+\cap\Phi_1$. Let $\De$ be the basis defined by $\Phi^+$ and put $\De_0=\De\cap\Phi_0$, $\De_1=\De\cap\Phi_1$. 

Denote by $\La\subset\lit^*$ the set of integral weights of $\Phi$ and observe that 
$\grs$ preserves $\grL$.
Let $\La^+$ be the set of dominant
weights with respect to $\Phi^+$ and 
let $\om_\al$ be the fundamental weight dual to the simple coroot $\coal$ for
$\al\in\De$.
For 
$\grl \in \grL^+$ let also $V_{\grl}$ 
be the irreducible representation of 
$\gog$ of highest weight $\grl$.

We say that $\grl\in \grL^+$ is \emph{spherical} if there exists 
$h \in V_\grl \senza \{0\}$ fixed by $\goh$
(i.e. $\goh \cdot h =0$): 
in this case the vector $h$ is also unique up to scalar and we denote it by $h_\grl$.
We denote the set of spherical weights by $\grO^+$
and we denote by $\grO$ the lattice generated by the spherical weights.

For a root $\gra$ define $\tal\doteq\al-\si(\al)$ and let $\wt \Phi=\{\tal \st \gra \in \Phi_1\}$.
This is a (not necessarily reduced) root system of rank $\ell$ with basis $\wt \grD = \{\wt \gra \st \gra \in \grD_1\}$. As proved in \cite{CM_SMT} using a result of Helgason, 
$\grO \cap \grL^+ = \grO^+$ and $\grO$  
can be identified with the lattice of integral weights of the root system 
$(\wt \Phi,\Omega \otimes_\mZ \mR)$. Given a weight $\grl$ we define the 
\emph{$\Om$ support} to be the set 
$\supp_\Om(\la) = \{\wt \gra \in \wt \grD \st (\la,\tal) \neq 0\}$.
We introduce also the lattice $\retrad$ generated by $\wt \Phi$ and the cone $\retrad^+ = \sum_{\wt \gra \in \wt \grD} \mN\,{\wt \gra}$.

Now we come to the construction of complete symmetric varieties following
De Concini and Procesi \cite{CP}.
Let $G$ be an algebraic group over $\mk$ whose Lie algebra is isomorphic to $\gog$.
The action of $\grs$ on $\gog$ lifts to an automorphism, still denoted by $\grs$,
of $G$. Let $H$ be the normalizer in $G$ of the Lie algebra $\goh \subset \gog$.
As explained in \cite{CP} $H$ is the maximal subgroup having $\goh$ as Lie algebra.
If $G$ is an adjoint group, $H$ coincides with
the fixed point set of $\grs$ in $G$. Hence $G/H$ is a symmetric variety.
However, since $G/H$ does not depend on the choice of the group $G$ over $\gog$ we will 
prefer to choose $G$ simply connected.
We introduce also the torus $T$ (resp. $T_0$ and $T_1$) whose Lie algebra is
$\got$ (resp. $\got_0$ and $\got_1$), and the parabolic subgroup $P$ of $G$ 
associated to $\grD_0$ (in general to a subset $I \subset \grD$ we associate the parabolic subgroup whose Lie algebra is given by $\got \oplus \bigoplus_{\gra \in \Phi_I\cup \Phi^+}\gog_\gra$, where $\Phi_I$ is the root subsystem of $\Phi$ generated by $I$).

In our study it will turn useful to consider also the degenerate case
$G=\{e\}$ or more generally $\grs=\id$. In this case, of course, $G/H$ is
just a single point.

Let now $\grl_1,\dots,\grl_m$ be spherical weights with disjoint
$\Om$ supports such that 
$\suppw(\grl_1)\cup \dots \cup \suppw(\grl_m)= \wt \grD$ and consider the point
$$
x_0 = ([h_{\grl_1}], \dots, [h_{\grl_m}]) \in \mP(V_{\grl_1})\times \dots \times \mP(V_{\grl_m}). 
$$
We define the variety $X=X(\grs)$ as 
$\overline{G x_0} \subset \mP(V_{\grl_1})\times \dots \times \mP(V_{\grl_m})$. 
Notice that $x_0$ is the unique point fixed by $H$ in $X$ and that the map $g \longmapsto gx_0$ 
induces an embedding $G/H \incluso X$ which is called the 
``minimal compactification'' of $G/H$.
Moreover the construction is independent on the choice of the weights $\grl_1,\dots,\grl_m$. 

We need also another description of the compactification.
Let $\grl$ be a spherical weight with $\suppw(\grl)=\wt \grD$ and 
consider a finite dimensional $\gog$ representation of the form  
$V\doteq V_{\grl} \oplus V'$.
Take $h=(h_{\grl}, h_{V'})\in V$ to be a vector fixed by $\goh$ such that all $T_1$ weights of $h_{V'}$ 
are of the form $\mu = \grl -\eta$ with $\eta \in \retrad^+$ 
and $\mu \neq \grl$.
Then, as proved in \cite{CP} \S 4, the map $G/H \ni gH \longmapsto g[h]\in \mP(V)$
extends to an isomorphism $X \lra \overline{G[h]}$.

The following Proposition describes the structure of the compactification.
\begin{proposition}[Theorem 3.1 in \cite{CP}]\label{pcompactification}
Let $X=X(\grs)$ be the compactification of $G/H$ described above then:
\begin{enumerate}[i)]
\item $X$ is a smooth projective $G$ variety;
\item $X\senza G\cdot x_0$ is a divisor 
with normal crossing and smooth irreducible components $S_1,\ldots,S_\ell$;
\item the $G$ orbits of $X$ correspond to the subsets of the indexes 
$1,2,\ldots,\ell$ so that the orbit closures are 
the intersections $S_{i_1}\cap S_{i_2}\cap\cdots\cap S_{i_k}$, with 
$1\leq i_1<i_2<\cdots<i_{k-1}<i_k\leq\ell$;
\item the unique closed orbit $Y\doteq\cap_{i=1}^\ell S_i$ is isomorphic to the partial flag variety $G/P$.
\end{enumerate}
\end{proposition}

We go on constructing some line bundles on the variety $X$.
Let $\la\in \grL^+$ be such that $\Pj(V_\la)$ contains a line $r$ fixed by $H$. 
One can show (see \cite{CP} and \cite{CS}) that the map
$G/H\ni gH\mapsto g\cdot r\in\Pj(V_\la)$ extends uniquely to a projection
$$
\psi_\la:X\rightarrow\Pj(V_\la).
$$
We define $\calL_\grl$ as the line bundle $\psi_\la^*\tlb$. 
If we restrict $\lb_\la$ on 
$G/P\simeq Y\hookrightarrow X$ we have the usual line bundle $G\times_P\field_{-\la}$ 
corresponding to $\la$ in the identification of 
$\Pic(G/P)$ with a sublattice of the weight lattice $\La$. Moreover we have
\begin{proposition}[Proposition 8.1 in \cite{CP}]\label{pinjectivepic} The map $\Pic(X)\rightarrow\Pic(Y)$ induced by the inclusion is injective.
\end{proposition}
So we can identify $\Pic(X)$ with a sublattice of the weight lattice. 
Further the line bundles constructed above account for all line bundles 
since we have
\begin{proposition}[Lemma 4.6 in \cite{CS}]\label{ppic1} $\Pic(X)$ corresponds to the lattice generated by the dominant weights 
$\la$ such that $\Pj(V_\la)^H$ is non void.
\end{proposition}
\noindent Notice that by construction every line bundle has a natural $G$ linearization.

Following the literature we introduce now a particular behavior
of a simple root. The action of the involution
$\si$ on the set of roots admits the following description. There exists an
involutive bijection $\osi:\De_1\rightarrow\De_1$ such that for every
$\al\in\De_1$ we have
$$
\si(\al)=-\osi(\al)-\be_\al
$$
where $\be_\al$ is a non negative linear combination of roots in $\De_0$.
We say that $\al\in\De_1$ is an \emph{exceptional} root if
$\osi(\al)\neq\al$ and $(\al,\si(\al))\neq0$. 
Notice that $\osi(\al)$ is exceptional if $\al$ is. Moreover
the compactification $X$ is said to be exceptional if
there exist exceptional roots.
\begin{proposition}[Theorem 4.8 in \cite{CS}]\label{ppic2} $\Pic(X)$ is generated by the spherical weights and the fundamental 
weights corresponding to the exceptional roots.
\end{proposition}

Now we describe the sections of a line bundle $\calL$ as a $G$-module. The first useful remark is that any irreducible $G$ module appears in $\Gamma (X,\calL)$ with multiplicity at most one (see Lemma~8.2 in \cite{CP}).

We analyze now the case of the divisors $S_i$. Let $\tal_1,\ldots,\tal_\ell$ 
be the elements of $\wt \grD$. 
Then, up to reindexing the $G$ stable divisors, we have
\begin{proposition}[Corollary 8.2 in \cite{CP}]\label{pdivisor} There exists a unique up to scalar $G$ invariant section $s_i\in \Gamma (X,\lb_{\tal_i})$ whose divisor is $S_i$.
\end{proposition}
For an element $\nu= \sum_{i=1}^{\ell} n_i {\wt \gra}_i \in \retrad^+$
the multiplication by $s^\nu \doteq \Pi_i s_i^{n_i}$ gives a linear map
$$
\Gamma (X,\lb_{\la-\nu})\rightarrow \Gamma (X,\lb_\la).
$$
If $\mu \in \Pic(X)$ is dominant then by construction of $\calL_\mu$ we certainly 
have a submodule of $\Gamma (X,\calL_\mu)$ isomorphic to $V_\mu^*$ obtained by the 
pull back of the homogeneous coordinates of $\mP(V_\mu)$ to $X$. Since, as already recalled, 
the multiplicity of
any irreducible submodule is at most one, we can speak of the submodule 
$V_\mu^*$ of $\Gamma (X,\calL_\mu)$ without ambiguity. 
If now $\grl \in \Pic(X)$ is any element such that $\grl-\mu \in \retrad^+$ we can 
consider the image of $V_\mu^*$ under the multiplication by $s^{\grl-\mu}$ from
$\Gamma (X,\calL_\mu)$ to $\Gamma (X,\calL_\grl)$. We call this image $s^{\grl -\mu}V_\mu ^*$.
We have the following Theorem:
\begin{proposition}[Theorem 5.10 in \cite{CP}]\label{pfiltration} Let $\la\in\Pic(X)$ then 
$$
\Gamma (X,\calL_\grl) = \bigoplus_{\mu \in (\grl-\retrad^+)\cap \grL^+}
s^{\grl -\mu}V_\mu ^*.
$$
Moreover for all $\nu = \sum_{i=1}^{\ell} n_i {\wt \gra}_i \in \retrad^+$ the 
set of sections vanishing on $S_i$ with multiplicity at least $n_i$ for 
$i=1,\dots,\ell$ is the image $s^{\nu}\Gamma (X,\calL_{\grl-\nu})$ of the multiplication 
by $s^{\nu}$ in $\Gamma (X,\calL_\grl)$.
\end{proposition}


\section{The stable subvarieties}

In this section we study the closures of $G$ orbits of $X$, that we call \emph{stable subvarieties}. Following De Concini and Procesi \cite{CP}, we review the structure of a stable subvariety, recalling, in particular, that such a variety is a fibration over a partial flag variety with fiber isomorphic to a complete symmetric variety. We use such result to prove Proposition \ref{prp:ind} lifting the surjectivity of the multiplication map from the fiber to a stable subvariety. This will be used in one inductive step for the proof of Theorem~A in next section.

We need to introduce some notation related with some special subgroups of $G$.
If $I$ is a subset of $\De$ containing $\grD_0$ and such that
$\bar \grs (I\cap \grD_1) \subset I\cap \grD_1$, let $\Phi_I \subset \Phi$ be the root subsystem of 
$\Phi$ generated by $I$ and define $G_I \subset G$ as the semisimple group 
associated to $\Phi_I$ (the subgroup whose Lie algebra is generated by $\gog_\gra$ 
for $\gra \in \Phi_I$). Also let $P_I = P \cap G_I$ be the parabolic subgroup of $G_I$ associated to $\De_0$
and let $T_I = T \cap G_I$ be a maximal torus of $G_I$. 
We call $\gog_I$ (resp. $\got_I$) the Lie algebra of $G_I$ (resp. $T_I$).
Observe that $\grs(\Phi_I)\subset \Phi_I$ as guaranteed by the assumption on $I$, hence
$G_I$ is stable under the action of $\grs$ and so we can define
$\grs_I \colon {G} _I \lra {G} _I$ as the restriction of $\grs$. 
We denote by $\goh_I$ the intersection $\goh \cap \gog_I$ and by $H_I$  
the normalizer of $\goh_I$ in $G_I$. Observe also that 
$\grs|_{\gog_\gra} = \id|_{\gog_\gra}$ for all $\gra \in \Phi_{\grD_0\cap I}$, hence $\dim (\got_I)_1$ is maximal.

We denote by $\grL_I$ the lattice of integral weights of $\Phi_I$ and 
by $\grL_I^+\subset\grL_I$ the set of dominant weights with respect to $I$.
Given two subsets $I,J$ of $\De$ such that $I\subset J$ we have a natural projection 
$r^J_I:\grL_J \lra \grL_I$ induced by the inclusion $\got_I\incluso \got_J$, 
and a canonical immersion $\mi^I_J : \grL_I \incluso \grL_J$ 
mapping a fundamental weight with respect to
$\Phi_I$ to the corresponding fundamental weight with respect to $\Phi_J$.
If we consider the case $J=\grD$ we see that 
$\grL_I$ is identified with the set of characters of $T_I$, 
hence also $G_I$ is simply connected. As recalled in the previous section, $H_J$ is the largest subgroup having $\goh_J$ as Lie algebra, then it is easy to see that $H_J \subset H$.

Now let $\wtI$ be a subset of $\tDe$ and set
$\ntI = \{ \gra \in \grD_1 \st \wt \gra \in \wtI \} \cup \grD_0$, notice that such $I$ satisfies the condition $\osi(I\cap\De_1)\subset I\cap\De_1$ above. 
Choose a one-parameter subgroup $\ga_{\wtI} \colon \field^* \lra T_1$ such that
$$
\bra \grg_{\wtI}, \wt\gra \ket = 0 \;\mif  \wt\gra \in \wt \grD \senza
\wtI \;\mand\;
\bra \grg_{\wtI}, \wt\gra \ket < 0 \;\mif  \wt\gra \in \wtI,
$$
where the pairing is given by the identification of a one parameter subgroup
of $T$ with an element of $\got$. We can now define the \emph{stable subvariety} $X_\wtI$ corresponding to $\wtI\subset\tDe$ as the closure $\wbar{G x_\wtI}$ of the orbit of the point $x_{\wtI} \doteq \lim_{t \lra 0} \grg_{\wtI}(t) x_0$.
\begin{prp}[Theorem 3.1 and Corollary 8.2 in \cite{CP}] We have
\begin{enumerate}[i)]
\item $X_{\{{\wt\gra}_i\}} = S_i$ for all ${\wt\gra}_i \in \wt \grD$;
\item $X_{\wtI \cup \wtJ} = X_{\wtI} \cap X_{\wtJ}$ and in particular $X_{\wt\grD}=Y$
is the unique closed $G$-orbit in $X$;
\item $X_{\wtI}$ is a projective smooth variety of dimension $\dim X - |\wtI|$.
\end{enumerate}
\end{prp}

In De Concini, Procesi \cite{CP} \S 5 the geometric structure of $X_{\wtI}$ 
is described. For the convenience of the reader we review below their results.

Fix $\wtI, \ntI$ as above and define $\wtJ = \wt \grD \senza \wtI$ and 
$\ntJ=(\grD \senza \ntI)\cup \grD_0$.
Let also $\retradJ^+ = \sum _{\wt \gra \in \wtJ} \mN \, \wt \gra$ and let $\retradJ$ be the 
lattice generated by $\retradJ^+$.
If $\grl \in \grL_J^+$ we denote by $W_{\grl}$ 
the irreducible representation of $G_J$ of highest weight $\grl$.

Let $\grl \in \grL$ be a spherical weight and consider a vector 
$h_{\grl}$ spanning the unique line in $V_\grl$
fixed by $H$. By \cite{CP} \S 2 we know that $h_\grl = v_\grl + \sum_{\mu \in \grl - (\retrad^+\senza \{0\})} w_{\mu}$ where $v_\grl$ is an highest weight vector and $w_\mu$ are eigenvectors of $T$ of weight $\mu$. We define
\begin{equation}\label{eq:hI}
h^{\wtI}_\grl \doteq v_\grl + \sum_{\mu \in \grl - (\retradJ^+ \senza \{0\})} w_\mu.\tag{$*$}
\end{equation}
Notice that, as a point of $X\subset \mP(V_{\grl_1})\times \dots \times \mP(V_{\grl_m})$  constructed in the previous section, we have
$$
x_{\wtI} = ([h_{\grl_1}^\wtI ] ,\dots,  [h_{\grl_m}^\wtI ])
$$
as one can easily see. 
\begin{lemma}\label{lhI}
If $\grl$ is a spherical weight then $[h_{\grl}^\wtI]$ is fixed by $H_J$. Further $h_{\grl}^\wtI = v_\grl$ if and only if $\suppw(\grl)\subset\wtI$.
\end{lemma}
\begin{proof}
Notice that $H_J \subset  H$ certainly fixes  $[h_\grl]$.
Observe also that $\grg_\calI$ commutes with $\gog_J$.
Indeed if $\gra \in \Phi_J$ and $e_\gra \in \gog_\gra$ then
$$
[\grg_\wtI , e_\gra ] = \bra \grg_\wtI , \gra \ket e_\gra =0.
$$
Hence $\grg_\calI$ commutes with $G_J$ and in particular with $H_J$.
The claim follows.

Suppose that $h^\wtI_\la=v_\la$. Let $\wtS = \wt \grD \senza \suppw(\grl)$ and 
$S=\{\gra \in \grD_1 \st \wt \gra \in \wtS \}\cup\grD_0$. Since 
$\mk h_{\grl}^\wtI$ is stable under $\goh_J$, we have 
$\goh_J\subset \stab_{\gog} \mk v_\grl = \got \oplus \bigoplus_{\gra \in \Phi_S \cup \Phi^+}\gog_\gra$. Now the inclusion $\suppw(\grl)\subset\wtI$ follows since, by the particular choice of $\Phi^+$, we have 
$\goh_J = \got_J \oplus \bigoplus _{\gra \in \Phi_0}\gog_\gra \oplus 
\bigoplus_{\gra \in \Phi_J\cap\Phi^+}\mk(e_\gra +\grs(e_\gra))$.

Now suppose that $\supp_\Om(\la)\subset\wtI$. 
Observe that if $\mu$ is a weight such that there exists $\gra \in I$
with $(\omega_\gra ,\grl-\mu) \neq 0$ then $\mu \notin \grl + \retrad_\calJ$.
Hence by formula \eqref{eq:hI} it is enough to 
show that for any weight $\mu\neq \grl$ with $(V_\grl)_\mu \neq 0$ 
there exists $\gra \in \grD$ such that $\bra \grl,\gra\cech \ket \neq 0$
and $(\omega_\gra ,\grl-\mu)>0$. 

We consider first the case $\mu = w\la\neq \grl$ with $w$ in the Weyl group of $\Phi$.
Let $L(\mu)$ be the minimum of the length of 
$w$ such that $\mu=w\la$. 
If $L(\mu)=1$ then $\mu = s_\gra (\grl) \neq \la$ for some $\gra$ and $\gra$ satisfies our
requests.

If $L(\mu)>1$ we proceed by induction: let $\mu = s_\gra \mu'$ 
where $\al$ is such that $L(\mu')=L(\mu)-1$. In particular 
$\grl-\mu = \grl-\mu ' + m \gra$ with $m>0$.
Then if $\grb \in \grD$ is such that 
$\bra \grb\cech,\grl \ket \neq 0$ and $(\omega_\grb,\grl-\mu')>0$ we have also
$(\omega_\grb,\grl-\mu)>0$.

Now in the general case we have that $\mu$ is in convex hull of $\{ w \grl\st w$ 
in the Weyl group of $\Phi\}$ and the claim follows.
\end{proof}
In order to describe $X_\wtI$ we consider a realization of $X$ as 
$\overline{G([h_\grl],[h_\mu])} \subset \mP(V_\grl) \times \mP(V_\mu)$
where $\grl$ and $\mu$ are two spherical weights such that $\supp_\Om(\grl)= \wtI$
and $\supp_\Om(\mu)= \wtJ$ (recall that $\wtJ=\tDe\senza\wtI$). 
Let $\pi : X \lra \mP(V_\grl)$ be the projection 
onto the first factor. By the previous Lemma \ref{lhI} we have $\pi (x_\wtI ) = [v_\grl ]$,
hence $\pi(X_\wtI) = G v_\grl \isocan G/Q$ where $Q$ is the parabolic 
subgroup of $G$ associated to $J$. We denote by $\pi_\wtI \colon X_\wtI \lra G/Q$ the restriction of $\pi$  to $X_\wtI$ and by $F_J$ the fiber of $\pi_\wtI$ over $[v_\grl]$.

\begin{proposition}[see \cite{CP} \S 5]\label{pfiber}
The fiber $F_J$ is the closure $\overline{{G}_J [h_\mu ^\wtI]}$. Moreover it is isomorphic to $X(\grs_J)$, the complete symmetric variety associated to $(G_J,\grs_J)$.
\end{proposition}
\begin{proof}
Consider the map $\pi_\wtI \colon X_\wtI \lra G / Q$. 
Since $\pi_\wtI$ is $G$ equivariant and $G/Q$ is homogeneous, $\pi_\wtI$ is a smooth map
and in particular the fiber $F_J$ is also smooth.
Observe that $F_J$ is closed under the action of $Q$, hence 
$x_{\calI'}=\lim _{t\to 0}\grg_{\calI'}x_\calI \in \overline {F_J \cap Gx_\calI }$ for all
$\calI' \supset \calI$.
Since $X_\calI = \bigcup_{\calI'\supset \calI}Gx_{\calI'}$, we find 
$F_J = \bigcup_{\calI'\supset \calI}(F_J \cap Gx_{\calI'}) = \overline{G x_\wtI \cap F_J}$. 
Hence $F_J = \overline{ Q x_\wtI }$ and $F_J$ is irreducible.
It is now easy to check that $\dim G_J/H_J = \dim X - |\wtI| - \dim G/Q=\dim F_J$. So 
${G}_J x_\wtI $ is dense in $F_J$ proving the first claim.

To prove the second claim we use the second construction of the compactification
$X(\grs_J)$ given in the previous section above Proposition \ref{pcompactification}.
Indeed observe that if $\nu \in \mu - (\retrad_\calJ^+ -\{0\})$ then $r(\nu)$ (the
restriction of the weight $\nu$ to $\got_J$) is strictly less then $r(\mu)$. Hence 
by \eqref{eq:hI} the vector $h_\mu^\calJ$ can be written in the form $h'_{r(\mu)}+h'_{W'}$
where $h'_{r(\mu)}$ is the only $\goh_J$ invariant vector of the $G_J$ module 
$W_{r(\mu)}$ and all $\got_J$ weights appearing in $h'_{W'}$ have weight strictly less 
than $r(\mu)$. 
\end{proof}
We apply this Proposition to the computation of $\Pic(X_\wtI)$.
As observed in the proof of the Proposition \ref{pfiber}, the stability of $F_J$ under the action
of $Q$ implies that $x_{\wt\grD} \in F_J$. Observe now that the action of $G_J$
on  $Y_J \doteq G_Jx_{\wt\grD}$ induces an identification of $G_J/P_J$ with $Y_J$ sending 
$P_J$ to $x_{\wt\grD}$, in particular $Y_J$ is the unique closed orbit of $F_J$. 
We identify $Y$ with $Gx_{\wt\grD}$ and we observe that the inclusion of 
$Y_J$ in $Y$ induces the natural inclusion of $G_J/P_J$ in $G/P$ sending $P_J=P\cap G_J$ to
$P$.
We denote by $\mj$ the inclusion of the closed orbit $Y\isocan G/P$ in 
$X_\wtI$ and with $\mj_J$ the inclusion of the closed orbit 
$Y_J\isocan G_J/P_J$ in $F_J$. We have the following commutative diagram (whose notation 
will be in force throughout the rest of this paper)
$$
\xymatrix{ 
0 \ar[r]& \Pic(G/Q)  \ar[dd]^\isocan \ar[r]^{\pi^*_\wtI} 
        & \Pic(X_\wtI) \ar[r]^{\mi_F^*}    \ar[d]^{\mj^*}       &\Pic(F_J) \ar[r]\ar@{^{(}->}[d]^{\mj_J^*} &0 \\
 & & \Pic(G/P)\ar[r]^{\mj_F^*}\ar[d]^\isocan & \Pic(G_J/P_J)\ar[d]^\isocan\\
0 \ar[r]& \grL_{\grD\senza J} \ar[r]^{\mi}
& \grL_{\grD\senza \grD_0}\ar[r]^{r} 
& \grL_{J\senza \grD_0 }\ar[r] & 0
}
$$
where:
\begin{enumerate}[i)]
\item $\mi = \mi^{\De\senza J}_{\De\senza\De_0}$, $r =r^{\De\senza\De_0}_{J\senza\De_0}$
and the third row is exact by construction.
\item $\mi_F$ is the inclusion of $F_J$ in $X_\wtI$. As in De Concini Procesi \cite{CP} there exists a one parameter subgroup of $G$ with only isolated
fixed points, hence $H^2(X_\wtI,\mZ) = \Pic(X_\wtI)$ and odd cohomology vanishes. So by the 
spectral sequence $H^p(G/Q, R^q\pi_{\wtI \, *} \mZ_{X_\wtI}) \then H^{p+q}(X_\wtI, \mZ_{X_\wtI})$ given by the fibration, we have that the first row is exact and $\Gamma (G/Q,R^2\pi_{\wtI \, *} \mZ_{X_\wtI})=H^2(F_J,\mZ)$ since $G/Q$ is simply connected.
\item $\mj_J^*$ is injective by Proposition \ref{pinjectivepic}.
\item $\mj_F^*$ is the pull back of the natural inclusion 
$G_J/P_J\hookrightarrow G/P$ that is compliant, as observed 
above, with the one induced by the inclusion $Y_J\subset Y$. Hence the square between the first and the second line is clearly commutative.
\item The isomorphisms mapping to the third row are the canonical 
identifications with the weight lattices.
\item The square on the left is commutative since $\mj\circ\pi_\wtI$ is 
the canonical projection induced by $P\subset Q$ and the pull back of 
the line bundle $G\times_{Q}\field_{-\la}$ gives the line bundle 
$G\times_{P}\field_{-\la}$. 
\item The square from the second and the third row is commutative since the 
line bundle $G\times_{P}\field_{-\la}$ restricted to $G_J/P_J$ 
gives the line bundle $G_J\times_{P_J}\field_{-\la |_{\got_J}}$.
\end{enumerate}
So $\mj^*$ is an injective map: if we identify $\Pic(G/Q)\simeq\La_{\De\setminus J}$ and $\Pic(X)$ with a sublattice of $\La_{\De\setminus\De_0}$ as in Section \ref{srecalls} then $\Pic(X_\wtI)$ is identified with the sublattice $\Pic(X)+\Pic(G/Q)$ of $\La_{\De\setminus\De_0}$.

Moreover observe that the inclusion
$\mi_X \colon X_\wtI \lra X$ induces an injective map $\mi_X^* \colon \Pic(X)\lra
\Pic(X_\wtI)$
and that, by the characterization of Proposition \ref{ppic2}, the map $\mi^*_F \circ\mi^*_X$ 
is surjective.
In particular we see that every line bundle on $X_\wtI$ has a natural 
$G$ linearization.
\begin{prp}\label{prp:sezXI} Let $\calL_\grl \in \Pic(X_\wtI)$. Then as a $G$-module we have 
$$
\Gamma (X_\wtI,\calL_\grl) = \bigoplus_ {\mu \in (\grl-\retradJ^+)\cap\grL^+}
s^{\grl-\mu}V_{\mu}^*.
$$
\end{prp}
\begin{proof}
Although,  for $\wtI \neq \vuoto$, this result it is not explicitly claimed in 
\cite{CP}, the proof of Theorem~8.3 there applies to this case without changes.
\end{proof}

We analyze now the relation between the sections $s_i$ of
the complete symmetric variety $X$ and the sections $s_{J,i}$ of the
complete symmetric variety $F_J$.

\begin{lem}\label{lem:sezsJ}
Up to rescaling the sections $s_{J,i}$ by a non zero constant factor we have 
$s_i|_{F_J} = s_{J,i}$  for all ${\wt \gra}_i \in \wtJ$.
\end{lem}
\begin{proof}
Observe that $s_i|_{F_J}$  and $s_{J,i}$ are $G_J$ invariant sections 
of the line bundle $\calL_{r({\wt \gra}_i)}$. Moreover $s_i|_{F_J}\neq 0$ 
by Proposition \ref{pdivisor} and the thesis follows. 
\end{proof} 

Looking at the decomposition in modules of the spaces of sections $\Gamma (X_\wtI,\lb_\la)$ and $\Gamma (F_J,\lb_{r(\la)})$ we see that they are indexed by the same weights. This suggests that one can prove the surjectivity of the multiplication map for sections on $X_\wtI$ using that on $F_J$. In order to make a rigorous proof out of this idea we analyze the lowest weight vectors as De Concini suggested us. 
Next lemmas prepare the work for this proof.

In the remaining of this section we will also make use of the following notation.
Set $U^-$ to be the unipotent subgroup of $G$ whose Lie algebra is 
$\bigoplus_{\gra\in \Phi^-}\gog_\gra$ and set $U_J^-=U^- \cap G_J$;
moreover if $\grl \in \grL^+$ let 
$V_\grl = \mk v_\grl \oplus V'_\grl$ be a $T$ stable 
decomposition of $V_\grl$ with $v_\grl$ an highest weight vector.

\begin{lem}\label{lem:ind1} If $\grl \in \grL^+_{\grD_1}$ then 
the restriction map
$$
\mj_F^* : \Gamma (G/P,\calL_\grl) \lra \Gamma (G_J/P_J,\calL_{r(\grl)})
$$
induces an isomorphism between $\Gamma (G/P,\calL_\grl)^{U^-}$ and 
$\Gamma (G_J/P_{J},\calL_{r(\grl)})^{U^-_J}$.
\end{lem}
\begin{proof}
Observe first that the map is ${G}_J$-equivariant, hence $\mj_F^*(\Gamma (G/P,\calL_\grl)^{U^-})\subset \Gamma (G_J/P_{J},\calL_{r(\grl)})^{U^-_J}$. 
Since they are both one dimensional vector spaces
it is enough to prove that $\mj_F^*\big(\Gamma (G/P,\calL_\grl) ^{U^-}\big) \neq \{0\}$.

The line bundle $\calL_\grl$ on $G/P$ can be constructed in the following way.
Consider the irreducible representation $V_\grl$ of highest weight 
$\grl$ and its highest weight vector $v_\la$. The stabilizer of $[v_\grl] \in \mP(V_\grl)$ 
contains $P$, hence 
we can construct a map $\psi \colon G/P \lra \mP(V_\grl)$ and 
$\calL_\grl = \psi^* \calO(1)$.
In particular $\Gamma (G/P,\calL_\grl) = V_\grl^*$ can be realized as the pull 
back through
the map $\psi$ of the space $V_\grl^*$ of coordinate functions on $\mP(V_\grl)$.

Let $\grf \in V^*_\grl$ be such that $\grf (v_\grl) =1$ and $\grf =0$ on $V'_\grl$,
then $\grf$ is a lowest weight vector in $V^*_\grl$, hence 
$\Gamma (G/P,\calL_\grl)^{U^-}=\mk \grf$. 
Observe that $\psi(\mj(P_J)) =\psi(P)=v_\grl$, hence 
$v_\grl \in \psi(\mj_F(G_J/P_J))$ and $\mj_F^*(\grf)\neq 0$.
\end{proof}

In the Lemma below we make use of the following straightforward consequence of the definition: the elements $r(\tal)$ with $\tal\in\wtJ$ form a basis of the restricted root system of $(G_J,\si_J)$.

\begin{lem}\label{lem:ind2} For all $\grl \in \Pic(X_\wtI)$ the restriction map 
$$
\mi_F^* : \Gamma (X_\wtI,\calL_\grl) \lra \Gamma (F_\ntJ,\calL_{r(\grl)})
$$
induces an isomorphism between $\Gamma (X_\wtI,\calL_\grl)^{U^-}$ and 
$\Gamma (F_\ntJ,\calL_{r(\grl)})^{U^-_J}$.
\end{lem}
\begin{proof} 
Observe first 
that the map $\mi_F^*$ is $G_J$ equivariant,
hence 
$\mi_F(\Gamma (X_\wtI,\calL_\grl)^{U^-}) \subset \Gamma (F_J,\calL_{r(\grl)})^{U^-_J}$.

Suppose now that $\grl$ is dominant and consider $V^*_\grl \subset \Gamma (X_\wtI,\calL_\grl)$ 
and $W_{r(\grl)} \subset \Gamma (F_\ntJ,\calL_{r(\grl)})$. 
Take $\grf \in V^*_\grl$ a lowest weight vector and observe that 
if $\mi_F^*(\grf)\neq 0$ then it is a vector in $\Gamma (F_\ntJ,\calL_{r(\grl)})$ of weight $-r(\grl)$, hence it spans $(W^*_{r(\grl)})^{U^-_J} \subset \Gamma (F_\ntJ,\calL_{r(\grl)})^{U^-_J}$. So it is enough to prove that $\grf|_{G_J/P_J}\neq 0$. Notice that, by the description 
of the sections of $\calL_\grl$ on $X_\wtI$ (Proposition \ref{pfiltration}), we have $\grf|_{G/P} \neq 0$ hence by the previous Lemma \ref{lem:ind1} we have $\grf|_{G_J/P_J} \neq 0$. In particular $\grf|_{F_J} \neq 0$.

Consider now the general case: let 
$M= (\grl - \retradJ^+) \cap \grL^+$
and 
$N= \big(r(\grl) - \sum_{\wt \gra \in \wtJ} \mN\,r(\tal)\big) \cap \grL^+_J$ 
and observe that there is a bijection 
between the two sets given by $M \ni \mu \longmapsto r(\mu) \in N$
and by 
$N \ni r(\grl) - \sum n_{\tal} r(\tal) 
\longmapsto \grl - \sum n_{\tal} {\wt\gra}\in M$. Notice that
$ \Gamma (X_\wtI,\calL_\grl)= \bigoplus_{\mu \in M} s^{\grl-\mu}V_\mu^*$, hence
$$
\Gamma (X_\wtI,\calL_\grl)^{U^-}= \bigoplus_{\mu \in M} \mk s^{\grl-\mu}\grf_\mu
$$
where $\grf_\mu \in V_\mu^*\subset \Gamma (X_\wtI,\calL_\mu) $ is a lowest weight vector.
Hence $\psi_{r(\mu)}\doteq\mi^*_F(\grf_\mu) \neq 0$ is a lowest weight vector in 
$W_{r(\mu)}^*\subset \Gamma (F_\ntJ,\calL_{r(\mu)})$ and, by Lemma \ref{lem:sezsJ}, 
$\mi^*_F(s^{\grl-\mu}\grf_\mu) =
s_\wtJ^{r(\grl-\mu)}\psi_{r(\mu)}$ up to a non zero scalar factor. 
Finally
\[
\mi^*_F(\Gamma (X_\wtI,\calL_\grl)^{U^-})= \bigoplus_{\mu \in M} \mk s_\wtJ^{r(\grl-\mu)}
\psi_{r(\mu)} = \bigoplus_{\nu\in N}\mk s_{\wtJ}^{r(\la)-\nu}\psi_\nu = \Gamma (F_\ntJ,\calL_{r(\grl)})^{U^-_J}.
\]
as claimed.
\end{proof}

If $\gol$ is a Lie algebra we denote by $\mbU(\gol)$ its universal enveloping
algebra and if $\gol$ is a subalgebra of $\gom$ then we consider
$\mbU(\gol)$ as a subalgebra of $\mbU(\gom)$. 
We introduce also the Lie algebra 
$\gou_J=\bigoplus_{\gra \in \Phi_J\cap\Phi^+}\gog_\gra$.

\begin{lem}\label{lem:ind3} Let $V, V'$ be two finite dimensional representations of $G$ and 
$W, W'$ be two finite dimensional representations of $G_J$.
Let $\phi\colon V \lra W$ (resp. $\phi' \colon V' \lra W'$) be a $G_J$ equivariant map
such that $\phi|_{V^{U^-}}$ (resp. $\phi'|_{{V'}^{U^-}}$) is an isomorphism between
$V^{U^-}$ and $W^{U_J^-}$ (resp. ${V'}^{U^-}$ and ${W'}^{U_J^-}$).
Then 
$$
\phi \otimes \phi' \big( (V\otimes V')^{U^-} \big) = (W\otimes W')^{U_J^-}.
$$
\end{lem}
\begin{proof}
Since the map is $G_J$ equivariant we certainly have that the left term is 
contained 
in the right one. Notice also that it is enough to study the case in which 
$V,V',W$ and 
$W'$ are irreducible.

Choose $\varphi$ and $\varphi'$ to be two lowest weight vectors of $V$ and $V'$ respectively.
Observe that, since $G_J$ is linearly reductive we can consider $W$ and $W'$ as the $G_J$ submodules generated by $\varphi$ and $\varphi'$ and the maps $\phi$ and $\phi'$ are just the projections along the $G_J$ invariant complements of $W$ and $W'$ in $V$ and $V'$. So we need to show that $(W\otimes W')^{U_J^-} \subset (V\otimes V')^{U^-}$.

Let $x \in (W\otimes W')^{U_J^-}$ and observe that it can be written as
$x = \sum_{m} e_m v \otimes e'_m v'$ where $e_m,e'_m \in \mbU(\gou_J)$.
To check  $x \in(V\otimes V')^{U^-}$  we verify $\gog_{-\gra} \cdot x = 0$ for
all $\gra \in\grD$. If $f_\gra \in \gog_{-\gra}$, we have
\begin{enumerate}[i)]
\item if $\gra \in \ntJ$ then $\gog_{-\gra} \subset \gog_J$ hence 
$f_\gra x = 0$ since $x \in (W\otimes W')^{U_J^-}$;
\item if $\gra \notin \ntJ$ then $f_\gra$ commutes with $e_m$ and $e'_m$, hence
$f_\gra y = \sum_m e_m f_\gra v \otimes e'_m v' +
\sum_m e_m v \otimes e'_m f_\gra v' = 0.$ \cvdqui
\end{enumerate}
\end{proof} 

We come now to the main result of this section.

\begin{prp}\label{prp:ind} Let $\grl,\mu \in \Pic(X_\wtI)$. If the multiplication map 
$$
\Gamma (F_\ntJ,\calL_{r(\grl)})\otimes \Gamma (F_\ntJ,\calL_{r(\mu)})\lra \Gamma (F_\ntJ,\calL_{r(\grl+\mu)})
$$
is surjective then also the multiplication map 
$$
\Gamma (X_\wtI,\calL_\grl)\otimes \Gamma (X_\wtI,\calL_\mu)\lra \Gamma (X_\wtI,\calL_{\grl+\mu})
$$
is surjective.
\end{prp}
\begin{proof} Consider the following commutative diagram where horizontal maps are given
by multiplication and vertical maps by restriction:
$$
\xymatrix{
\Gamma (X_\wtI,\calL_\grl)\otimes \Gamma (X_\wtI,\calL_\mu) \ar[r] \ar[d]   &
\Gamma (X_\wtI,\calL_{\grl+\mu}) \ar[d] \\
\Gamma (F_\ntJ,\calL_{r(\grl)})\otimes \Gamma (F_\ntJ,\calL_{r(\mu)}) \ar[r] &
\Gamma (F_\ntJ,\calL_{r(\grl+\mu)}).
}
$$
If we look at $U^-$ and $U^-_J$ invariants we obtain:
$$
\xymatrix{
\big(\Gamma (X_\wtI,\calL_\grl)\otimes \Gamma (X_\wtI,\calL_\mu)\big)^{U^-}  \ar[r] \ar@{->>}[d]   &
\Gamma (X_\wtI,\calL_{\grl+\mu})^{U^-}  \ar[d]^{\isocan} \\
\big(\Gamma (F_\ntJ,\calL_{r(\grl)})\otimes \Gamma (F_\ntJ,\calL_{r(\mu)})\big) ^{U_J^-} \ar@{->>}[r] &
\Gamma (F_\ntJ,\calL_{r(\grl+\mu)})^{U_J^-} 
}
$$
where, by Lemma \ref{lem:ind2}, the vertical map on the right is an isomorphism, 
by Lemma \ref{lem:ind2} and Lemma \ref{lem:ind3},
the vertical map on the left
is surjective 
and, by 
the assumption on the multiplication on $F_\ntJ$  
and the fact that $G_J$ is a linearly reductive group,
 the horizontal map on the bottom is surjective.

Hence also the horizontal map on the top has to be surjective. Further, since $\Gamma (X_\wtI,\calL_{\grl+\mu})^{U^-}$
generates $\Gamma (X_\wtI,\calL_{\grl+\mu})$ as a $G$-module we deduce that the multiplication
map on $X_\wtI$ is surjective too.
\end{proof}


\section{Projective normality}\label{sprojectivenormality}

In this section we prove the surjectivity of the multiplication map for dominant line bundles, i.e. for line bundles $\lb_\la$ with $\la\in\Pic^+(X)\doteq\Pic(X)\cap\La^+$. The main ingredients will be Proposition \ref{prp:ind} that allow us to set up an induction on the dimension of the symmetric variety, and Theorem~B (proved in Section \ref{slowtriple}).

We denote by $W$ the Weyl group of $\Phi$ and by $w_0$  the longest element of $W$.
Given two weights $\la,\mu\in\Pic(X)$, we call $m_{\la,\mu}$ the multiplication map from $\Gamma (X,\calL_{\grl})\otimes \Gamma (X,\calL_{\mu})$ to $\Gamma (X,\calL_{\grl+\mu})$. The following lemma deals with a very special case of Theorem~A.
\begin{lemma}\label{linvariant}
Given a weight $\la\in\Pic^+(X)$ set $\mu\doteq-w_0\la$. Then $\mu\in\Pic^+(X)$ and $s^{\grl+\mu}V_0^*\subset\Im m_{\la,\mu}$.
\end{lemma}
\begin{proof}
Observe first that we can identify $V_{\grl}^*$ with $V_{\mu}$. By Proposition \ref{ppic1} there exists $h_\la\in V_\la\setminus\{0\}$ such that $[h_\la]\in\Pj(V_\la)^H$. Hence there exists a one dimensional $H$-submodule $\chi$ of $V_\la$. Being $H$ reductive we have that $\chi^*$ is a submodule of $V_\mu$. So there exists 
$h_\mu\in V_\mu\setminus\{0\}$ such that $[h_\mu]\in\Pj(V_\mu)^H$ and $\bra h_{\grl},h_{\mu} \ket =1$. By Proposition \ref{ppic1}, we deduce that $\mu\in\Pic(X)$. Clearly $\mu\in\La^+$, hence $\mu\in\Pic^+(X)$.

Now complete the vectors 
$h_\mu\in V_\la^*$ and $h_\la\in V_\mu^*$ to dual bases $h_\mu,v_1,\ldots,v_n$ and $h_\la,w_1,\ldots,w_n$. Consider the following element of $V_\la^*\otimes V_\mu^*$: 
$$
F=h_\mu\otimes h_\la+\sum_{i=1}^n v_i\otimes w_i.
$$
If we identify $V_\la^*\otimes V_\mu^*$ with $\End (V_{\grl})$, the element $F$ corresponds to the identity map, in particular it is a $G$-invariant vector. Hence $f=m_{\la,\mu}(F)\in \Gamma (X,\lb_{\la+\mu})$ is 
$G$-invariant. 
We claim that $f\neq0$ proving the lemma. Indeed consider the injection $\psi:X\hookrightarrow\Pj(V_\la)\times\Pj(V_\mu)$ defined by $\psi(x)=(\psi_\la(x),\psi_\mu(x))$ (see Section \ref{srecalls} for the definition of $\psi_\la,\psi_\mu$) and notice that $([h_\la],[h_\mu])\in\Im\psi$. We have
$$
f([h_\la],[h_\mu])=h_\mu([h_\la])h_\la([h_\mu])+\sum_{i=1}^n v_i([h_\la])w_i([h_\mu])=1
$$
since we have chosen dual bases.
\end{proof}
Recall that a complete symmetric variety $X=X(\si)$ is said to be \emph{simple} if $\lig$ has no $\si$ 
stable proper ideal. It is known (see for example Table~VI cap.~X in \cite{H}) that a simple complete symmetric variety corresponds to an irreducible root system $\tPhi$ and either $G$ is simple or $X$ is the compactification of a simple group. Further any complete symmetric variety is the
product of simple complete symmetric varieties.

We need to make some preliminary remark on the various lattices and on the relations between 
the Weyl group $W$ of $\Phi$ and the Weyl group $\wt W$ of $\tPhi$.
As claimed in Section \ref{srecalls}, the lattice $\grO$ generated by spherical
weights is identified with the set of integral weights of the reduced root system $\tPhi$. 
Further the set of spherical weight $\grO^+$ is equal to $\grL^+\cap \grO$ and corresponds
to the dominant chamber defined by $\wt \grD$.
So the longest element $\wt w _0$ of the Weyl group $\wt W$ and
the longest element $w _0$ of the Weyl group $W$ act in the same way on $\grO$.

Before giving the proof of the surjectivity of $m_{\grl,\mu}$ we introduce some 
notation and make some remarks to treat the exceptional case. 
If $X$ is an exceptional simple complete symmetric variety then $\tPhi$ is of type $\typeBC_\ell$ 
and there exist exactly two (simple) exceptional roots that we call $\gra$ and $\grb$.
Also we denote by $\tal_\ell$ the unique simple root in $\wt \grD$ such that $2\tal_\ell
\in \tPhi$ and by $\tom_\ell$ the fundamental weight dual to $(2\tal_\ell)\cech$ (see part
2 in Section \ref{slowtriple}). We have $\Pic(X) =\grO \oplus \mZ \om_\al$ and 
$\Pic^+(X) =\grO^+ + \mN\,\om_\al + \mN\,\om_\grb$. Moreover $\om_\gra +\om_\grb=\tom_\ell$ and 
$-w_0 (\om_\al) =\om_\grb $ (see \cite{CM_SMT}).
We remind also that if $X$ is non exceptional then $\Pic(X)=\grO$ and $\Pic^+(X)=\grO^+$.

Finally we notice that, given two weights $\grl,\mu\in\Pic(X)$, we have $\mu \in (\grl -\retrad^+)$ if and only if $\mu \leq \grl$ with respect to the dominant order of $\tPhi$ (see part~3 in Section \ref{slowtriple} for the definitions in the exceptional case).

Now we come to the main result of this paper.
\begin{theoremA}\label{tsurjectivity}
Let $\la,\mu$ be two weights in $\Pic^+(X)$. Then the multiplication map
$$
m_{\la,\mu}:\Gamma (X,\lb_\la)\otimes \Gamma (X,\lb_\mu)\rightarrow \Gamma (X,\lb_{\la+\mu})
$$
is surjective.
\end{theoremA}
\begin{proof}
We prove first the case in which $\grl, \mu \in \grO^+$.
We proceed by induction on $\dim X$ and on the dominant order on $\la$ and $\mu$ with respect to the reduced root system $\tPhi$. If $\dim X=0$, i.e. $X$ is a point, then  $\la=\mu=0$ and the claim is obvious. Also if $X$ is not simple, say $X=X_1\times X_2$ with $\dim X_1,\dim X_2>0$ then the thesis follows by induction on the dimension using the description of the sections of the line bundles on $X$ in Proposition \ref{pfiltration}. So we can assume that $X$ is simple, hence that $\tPhi$ is irreducible.

We fix a notation. Given a weight $\eta\in\Pic^+(X)$ set $\La(\eta)=\La^+\cap(\eta-\wt R^+)$. Notice that clearly $\Gamma (X,\lb_\eta)=\oplus s^{\eta-\nu}V_\nu^*$ where the sum runs over $\nu\in\La(\eta)$. Our thesis is that $s^{\la+\mu-\nu}V_\nu^*\subset\Im m_{\la,\mu}$ for all $\nu\in\La(\la+\mu)$. We divide the set of $\nu\in\La(\la+\mu)$ in three different classes.

\emph{First class. } This is the class of weights $\nu\in\La(\la+\mu)$ such that the following condition is verified: there exists $\la',\mu'\in\Om^+$ such that $(\la',\mu')\neq(\la,\mu)$, $\la'\in\La(\la)$, $\mu'\in\La(\mu)$ and $\nu\in\La(\la'+\mu')$. Consider the following commutative diagram
$$
\xymatrix{
\Gamma (X,\lb_{\la'})\times \Gamma (X,\lb_{\mu'})\ar[r]^(.6){m_{\la',\mu'}}\ar[d]^{s^{\la-\la'}\otimes s^{\mu-\mu'}} & \Gamma (X,\lb_{\la'+\mu'})\ar[d]^{s^{\la-\la'+\mu-\mu'}}\\
\Gamma (X,\lb_\la)\times \Gamma (X,\lb_\mu)\ar[r]^(.6){m_{\la,\mu}} & \Gamma (X,\lb_{\la+\mu}).\\
}
$$
Notice that $m_{\la',\mu'}$ is surjective by induction on the dominant order on $\la$ and $\mu$. Also, $s^{\la+\mu-\nu}V_\nu^*$ is contained in the image of the right vertical map since $s^{\la'+\mu'-\nu}V_\nu^*\subset \Gamma (X,\lb_{\la'+\mu'})$. So $s^{\la+\mu-\nu}V_\nu^*$ is contained in $\Im m_{\la,\mu}$.

\emph{Second class.} This class is formed by the weights $\nu\in\La(\la+\mu)$ such that $\la+\mu-\nu=\sum_{i=1}^\ell c_i\tal_i$ with an index $i$ such that $c_i=0$.
 
If we consider the restriction of sections to  
the stable subvariety
$X_{\{\tal_i\}}$ we have the following commutative diagram:
$$
\xymatrix{
\Gamma (X,\lb_\la)\times \Gamma (X,\lb_\mu)\ar[r]^(.6){m_{\la,\mu}}\ar[d] & \Gamma (X,\lb_{\la+\mu})\ar[d]\\
\Gamma (X_{\{\tal_i\}},\lb_\la)\times \Gamma (X_{\{\tal_i\}},\lb_\mu)\ar[r]^(.6){m_{\la,\mu}} & \Gamma (X_{\{\tal_i\}},\lb_{\la+\mu}).\\
}
$$
Notice that the bottom horizontal map is surjective by Proposition \ref{prp:ind} and induction on dimension since $\dim(F_J)<\dim X$, where $J=\{\al\in\De_1\ |\ \tal\neq\tal_i \}\cup\De_0$. Also notice that the left vertical map is surjective by Lemma \ref{lem:ind2}. 
Finally observe that $s^{\la+\mu-\nu} V_\nu^*$ appears in the decomposition of $\Gamma (X_{\{\tal_i\}},\lb_{\la+\mu})$ described in Proposition \ref{prp:sezXI} since $c_i=0$. 
Hence $s^{\la+\mu-\nu}V_\nu^*\subset\Im m_{\la,\mu}$ since $G$ is reductive and all modules appears with multiplicity at most one.

\emph{Third class.} This is the set of the remaining $\nu\in\La(\la+\mu)$: what is left is described by the triples of weights $(\la,\mu,\nu)$ of $\Om^+$ such that i) $\la'\in\La(\la)$, $\mu'\in\La(\mu)$, $\nu\in\La(\la'+\mu')$ implies $\la'=\la$ and $\mu'=\mu$ and ii) $\la+\mu-\nu=\sum_{i=1}^n c_i\tal_i$ with $c_i\geq1$ for all $i=1,\ldots\ell$. So $(\la,\mu,\nu)$ is a low triple for the root system $\tPhi$.

By Theorem~B applied to the root system $\tPhi$ we have that $\nu=0$ and that $\mu= - w_0 \grl$. 
So we can conclude using Lemma \ref{linvariant}. 

This concludes the proof for $\grl,\mu\in \grO^+$ and in particular the case of $X$ non exceptional. 
So we are reduced to study the case of $X$ exceptional and, using what already done, we can assume the theorem true for all $\grl, \mu \in \grO^+$.
We proceed again by induction on $\dim X$ and on the dominant order on $\Pic^+(X)$ defined
by $\tPhi$. As in the case of $\grO^+$ we can assume that $X$ is a simple
exceptional complete symmetric variety. 

We consider first a particular case. Let $\grl = \om_\gra$ and $\mu=\mu' + h\om_\al$ for some $h\geq 0$ and $\mu'\in\Om^+$. By Proposition \ref{pfiltration} and Lemma \ref{lenlargedbelow} we have the following decompositions:
\begin{align*}
\Gamma (X,\lb_{\om_\al})    &=V_{\om_\al}^*,\\
\Gamma (X,\lb_\mu)          &=\oplus_{\nu \in \grL(\mu)} s^{\mu-\nu} V_{\nu}^*,\\
\Gamma (X,\lb_{\om_\al+\mu})&=\oplus_{\nu \in \grL(\mu)} s^{\mu+\om_\gra -\nu}V_{\om_\gra+\nu}^*.
\end{align*}
Denote by $\grf_\eta$ a highest weight vector of the module $V_\eta^*$. So to prove the surjectivity of $m_{\om_\al,\mu}$ in this situation, we observe
that $\grf_{\om_\al}\otimes 
s^{\mu-\nu}\grf_{\nu}$ is a non zero highest weight vector of weight 
$\om_\gra+\mu$ hence a multiple of $s^{\mu-\nu}\grf_{\nu}$. 

Now we proceed as in the case of $\grO^+$ and we observe that
the arguments given in the ``first class" and in the ``second class" above holds without any change and so we are reduced again to study the low triples of weights $(\la,\mu,\nu)$ of $\Pic^+(X)$ with respect to
the dominant order defined by $\tPhi$.
Hence, by Proposition \ref{pexclowtriple},  we have that, up to symmetry, 
$\grl = a \om_\gra$ and $\mu = b \om_\grb$
for some integers $a\geq b \geq 1$.

We analyze first the particular case $a=b=1$. In this case, by Proposition \ref{pexclowtriple},
we have also $\nu=0$. 
Hence we can apply Lemma \ref{linvariant} since $\om_\be=-w_0\om_\al$ and conclude that 
$s^{\tom_\ell}V_0^*$ is contained in the image of $m_{\om_\al,\om_\be}$ as claimed. 

Now we deduce the surjectivity in the remaining cases of $\grl=a\om_\gra$ and $b=b\om_\grb$ with 
$a\geq b\geq 1$ using the associativity of multiplication.
We consider the following commutative diagram 
$${\small
\xymatrix{
\Gamma (X,\lb_{\om_\al})^{\otimes (a-b)}\otimes 
\big( \Gamma (X,\lb_{\om_\al}) \otimes \Gamma (X,\lb_{\om_\be})\big)^{\otimes b}
\ar@{=}[r] \ar@{->>}[d]_{m_1} 
& \Gamma (X,\lb_{\om_{\al}})^{\otimes a}\otimes \Gamma (X,\lb_{\om_{\be}})^{\otimes b}
\ar[d]  \\
\Gamma (X,\lb_{\om_\al})^{\otimes(a-b)}\otimes (\Gamma (X,\lb_{\tom_\ell}))^{\otimes b}
\ar@{->>}[d]_{m_2}  & 
\Gamma (X,\lb_{a\om_{\al}})\otimes \Gamma (X,\lb_{b\om_{\be}}) \ar[d]_{m_{a\om_\gra,b\om_\grb}} \\
\Gamma (X,\lb_{\om_\al})^{\otimes(a-b)}\otimes \Gamma (X,\lb_{b\tom_\ell})
\ar@{->>}[r]^(.55){m_3}
& \Gamma (X,\lb_{(a-b)\om_\al+b\tom_\ell}).
}}
$$
where 
\begin{enumerate}[i)]
\item $m_1 = \id \otimes m_{\om_\gra,\om_\grb}^{\otimes b}$ is surjective using what proved in the particular case $a=b=1$;
\item $m_2$ is the multiplication of the sections in $\Gamma (X,\calL_{\tom_\ell})^{\otimes b}$ that 
is surjective 
using what proved in the particular case of $\grl, \mu \in \grO^+$;
\item $m_3$ is the multiplication of sections in $\Gamma (X,\calL_{\om_\gra})^{\otimes (a-b)}$ and
$\Gamma (X,\calL_{b\tom_\ell})$ that is surjective using what proved in the case $\grl = \om_\gra$
and $\mu=\mu' +h\om_\gra$ analyzed above.
\end{enumerate}
Hence also $m_{a\om_\al,b\om_\be}$ is surjective. 
\end{proof}

We see an example showing the necessity to assume the line bundles to be dominant for the surjectivity of the multiplication map. Simply take a complete symmetric variety with $\tPhi$ of type $\typeA_2$ and let $\la=\tal_1$, $\mu=\tal_2$. Then $\Im m_{\la,\mu}=s_1s_2V_0^*$, whereas $\Gamma (X,\lb_{\la+\mu})=V^*_{\tom_1+\tom_2}\oplus s_1s_2V^*_0$. As another example with $\la=\mu$, consider the case of the compactification of the group of type $\typeC_\ell$ and let $\la=\mu=-\om_{\ell-1}+\om_\ell$. Then $\Gamma (X,\lb_\la)=\Gamma (X,\lb_\mu)=0$ while $\Gamma (X,\lb_{\la+\mu})=s_\ell V_0^*$ as one can easily see by Proposition \ref{pfiltration}.

\begin{cor} For all $\calI \subset \wt \grD$ and for all $\grl,\mu \in \Pic^+(X_\calI)= \grL^+ \cap \Pic(X_\calI)$ we have that the multiplication map 
$$
\Gamma (X_\calI,\lb_\la)\otimes \Gamma (X_\calI,\lb_\mu)\rightarrow \Gamma (X_\calI,\lb_{\la+\mu})
$$
is surjective.
\end{cor}

\begin{proof} This follows at once by Theorem~A and Proposition \ref{prp:ind}.
\end{proof}

As a consequence (see for example Hartshorne Exercise II.5.14) we have
\begin{corollary}\label{cnormality}
Let $\calI$ be a subset of $\tDe$ and let $\lb_\la\in\Pic^+(X_\calI)$ be a dominant line bundle. Consider the map $X_\calI\rightarrow\mP(\Gamma (X_\calI, \calL_{\grl})^*)$ defined by the line bundle $\lb_\la$. Then  the cone over the image of $X_\calI$ is normal. In particular this apply to $X$.
\end{corollary}


\section{\Low triples}\label{slowtriple}

In this section we introduce and study \low triples for an irreducible root system. As seen in the previous section, these are the triples of dominant weights that furnish the base step for the inductive proof of the surjectivity of the multiplication map. In the first part we consider only reduced root systems, developing a bit of general theory as far as we are able to. Then we consider the case of a non reduced root system that present no more real difficulty using what already done in the first part. Finally we give a little more general and technical result concerning \low triples for an enlarged weight lattice for type $\typeBC$; this is needed for exceptional complete symmetric varieties.

\subsection{Reduced root system}\mbox{}\\
Let $\Phi$ be an irreducible reduced root system with positive roots $\Phi^+$, relative base $\De=\{\al_1,\ldots,\al_\ell\}$, and let $\La$ be the set of integral weights for $\Phi$. We denote by $\om_1,\ldots,\om_\ell$ the base of $\La$ dual to the coroots $\al_1^\vee,\ldots,\al_\ell^\vee$ and by $\La^+$ the $\N$ cone of dominant weights. Let $\leq$ be the usual dominant order on $\La$, i.e. the order: $\mu\leq\la$ if and only if $\la-\mu\in R^+$ where $R^+\subset R$ is the $\N$ cone generated by the simple roots and $R$ is the lattice generated by the simple roots. If $w_0$ is the longest element of the Weyl group of $\Phi$, we say that two dominant weights $\la$, $\mu$ are \emph{dual} to each other if $\mu=-w_0\la$. As a notation, in the sequel we use the numbering of the simple roots in \cite{B}.

We denote by $[\ell]$ the set of indexes $\{1,\ldots,\ell\}$. Given a weight $\la=\sum_{i\in[\ell]}a_i\om_i$, we define the \emph{support} as the set $\supp(\la)$ of indexes $i$ such that $a_i\neq0$ and $\supp^+(\la)$ as the set of indexes $i$ such that $a_i>0$. We define also the \emph{height} of the weight $\la$ as $\height(\la)=\sum_{i\in[\ell]}a_i$. Let $I_0$ be the maximal simply laced subsystem of $[\ell]$ containing $1$, and denote by $\hts(\la)$ either $\sum_{i\in I_0, a_i>0}a_i$ if $\Phi$ is not of type $\typeF_4$, or the maximum of $\sum_{i\in\{1,2\}, a_i>0}a_i$ and $\sum_{i\in\{3,4\}, a_i>0}a_i$ otherwise. Given an element of the root lattice $\be=\sum_{i\in[\ell]}b_i\al_i$, we define $\supp_\Phi(\be)$ as the set of indexes $i$ such that $b_i\neq0$. Finally let $\zeta\doteq\al_1+\al_2+\cdots+\al_\ell$.

Now we see the definition of the main object of this section
\begin{definition}\label{dlowtriple}
Let $\la,\mu,\nu$ be three dominant weights. We call $(\la,\mu,\nu)$ a \emph{\low triple} if the following two conditions hold:\\
i) if $\la',\mu'$ are dominant weights such that $\la'\leq\la$, $\mu'\leq\mu$ and $\nu\leq\la'+\mu'$ then $\la'=\la$ and $\mu'=\mu$,\\
ii) $\nu+\sum_{i=1}^\ell\al_i\leq\la+\mu$. 
\end{definition}
Notice that \low triples are of a very special kind. Indeed $\nu$ must be ``very close'' to $\la+\mu$ by i), whereas $\nu$ must be ``quite far'' from $\la+\mu$ by ii). The rest of this section is devoted to give precise meaning to this idea proving the following
\begin{theoremB}
$(\la,\mu,\nu)$ is a \low triple if and only if $\la$ and $\mu$ are minuscule weights dual to each other and $\nu=0$.
\end{theoremB}
For $I\subset[\ell]$, let $\Phi_I$ denote the root subsystem generated by the simple root $\al_i$ for $i\in I$. If $\Phi_I$ is irreducible we denote by $\theta_I$ the highest short root (we notice that in \cite{Stem} the root $\theta_I$ is called the \emph{local short dominant root} relative to $I$), considering all roots short if $\Phi_I$ has just one root length. In particular $\theta=\theta_{[\ell]}$ is the highest short root of $\Phi$.

Given two weights $\la,\mu\in\La^+$, we say that $\la$ \emph{covers} $\mu$ if i) $\mu\leq\la$ and ii) $\mu\leq\eta\leq\la$, $\eta\in\La^+$ imply $\eta=\mu$ or $\eta=\la$. Of particular importance for the sequel is the following characterization of the covering relation for the dominant order $\leq$ restricted to the set of dominant weights (see \cite{Stem} for details).
\begin{proposition}[Theorem~2.6 in \cite{Stem}]\label{pcovering}
If $\mu$ covers $\la$ in $(\La^+,\leq)$ and $I=\supp_\Phi(\mu-\la)$, then $\Phi_I$ is irreducible and either $\mu-\la=\theta_I$, or $\Phi_I=\Phi\simeq\typeG_2$ and $\mu-\la=\al_1+\al_2$.
\end{proposition}
We explicitly write down the highest short root for each type of root system. This is an easy computation using the tables in \cite{B}.
\begin{center}
\begin{tabular}{c|c}
type of $\Phi$ & highest short root\\
\hline
$\typeA_\ell$ & $\al_1+\cdots+\al_\ell=\om_1+\om_\ell$\\
$\typeB_\ell$ & $\al_1+\cdots+\al_\ell=\om_1$\\
$\typeC_\ell$ & $\al_1+2(\al_2+\cdots+\al_{\ell-1})+\al_\ell=\om_2$\\
$\typeD_\ell$ & $\al_1+2(\al_2+\cdots+\al_{\ell-2})+\al_{\ell-1}+\al_\ell=\om_2$\\
$\typeE_6$ & $\al_1+2\al_2+2\al_3+3\al_4+2\al_5+\al_6=\om_2$\\
$\typeE_7$ & $2\al_1+2\al_2+3\al_3+4\al_4+3\al_5+2\al_6+\al_7=\om_1$\\
$\typeE_8$ & $2\al_1+3\al_2+4\al_3+6\al_4+5\al_5+4\al_6+3\al_7+2\al_8=\om_8$\\
$\typeF_4$ & $\al_1+2\al_2+3\al_3+2\al_4=\om_4$\\
$\typeG_2$ & $2\al_1+\al_2=\om_1$\\
\end{tabular}
\end{center}

Given a subset $I\subset[\ell]$ we say that $I$ is irreducible (or of type $\typeA$, $\typeB$, $\ldots$) if $\Phi_I$ is.
We make the following remark reading it out of the table. If $I\subset[\ell]$ is irreducible and not of type $\typeA$, then $\theta_I=\om_i-\eta$ for some $i\in[\ell]$ and $\eta\in\La^+$ with $I\cap\supp(\eta)=\vuoto$. So we can define a map $f$ from the irreducible subsets of $[\ell]$ of type not $\typeA$ to $[\ell]$ mapping $I$ to the index $f(I)=i$ above.

Another direct consequence of the computation of highest short roots is
\begin{lemma}\label{lhsrcomparison}
Let $I,J\subset[\ell]$ not of type $\tpA$. 
Then $J\subset I\subset[\ell]$ if and only if $\theta_I\geq\theta_J$.
\end{lemma}

We need also the following criterion.

\begin{lem}\label{luseful}
Let $(\grl,\mu,\nu)$ be a triple of dominant weights and suppose that there exists 
$\grg=\sum_{i\in[\ell]} c_i \gra_i$ with $c_1,\dots,c_\ell$ non negative integers, such that 
\begin{enumerate}[i)]
\item $\grl - \grg \in \grL^+$; 
\item $\grl+\mu - \nu \geq \grg$.
\end{enumerate}
Then $(\la,\mu,\nu)$ is not a \low triple.
\end{lem}
\begin{proof}
Set $\la'\doteq\la-\grg$. 
Then $\nu\leq\la'+\mu$ against the first condition in Definition \ref{dlowtriple} of \low triple.
\end{proof}

In the sequel we will use Lemma \ref{luseful} above in the case $\sum_{i\in[\ell]}c_i\al_i=\theta_I$ 
for some irreducible subset $I\subset[\ell]$, and in the case in which there exists $I\subset[\ell]$ 
such that $\la-\sum_{i\in I}\al_i\in\La^+$. Indeed notice that in this latter case hypothesis ii) in 
the lemma is guaranteed by $\la+\mu-\nu\geq\zeta\geq\sum_{i\in I}\al_i$.

We begin the proof of Theorem~B with the following lemma limiting the possibilities for \low triples.
\begin{lemma}\label{lheight}
Suppose that $(\la,\mu,\nu)$ is a low triple, let $\be\doteq\la+\mu-\nu$. Then we have $\height_0^+(\be)\leq 2$.
\end{lemma}
\begin{proof}
First we treat the case of $\Phi$ not of type $\typeF_4$. Suppose that $\height_0^+(\be)\geq3$. 
Then $\hts(\la)+\hts(\mu)=\hts(\la+\mu)\geq\hts(\be)\geq3$. 
Hence either $\hts(\la)\geq2$ or $\hts(\mu)\geq2$. By symmetry we can suppose $\hts(\la)\geq2$.

We have two possibilities 1) $\la=2\om_j+\eta$ for some $j\in[\ell]$ and $\eta\in\La^+$ or 2) $\la=\om_i+\om_j+\eta$ for some $i,j\in I_0$ and $\eta\in\La^+$. In the first case $\la-\al_j\in\La^+$ and $\be\geq\zeta=\sum_{i\in[\ell]}\al_i>\al_j$, and this is impossible $(\la,\mu,\nu)$ being a low triple.

In the second case let $J\subset I_0$ be the (type $\typeA$) segment connecting $i$ and $j$ in the Dynkin diagram of $\Phi$. Then $\la-\theta_J=\la-\om_i-\om_j+\eta'$ for some $\eta'\in\La^+$. Hence $\la-\theta_J\in\La^+$ and $\be\geq\theta_J=\sum_{h\in J}\al_h$. This is impossible.

Finally if $\Phi$ is of type $\typeF_4$ we argue as above taking into account the definition of $\hts$ for such type.
\end{proof}

Given an element $\be$ of the root lattice $R$ we denote by $\calI(\be)$ the set of $I\subset[\ell]$ such that $I$ is irreducible of type not $\typeA$ and $\be\geq\theta_I$.
We can consider such set ordered by the inclusion relation, notice that this is also the order on the corresponding highest short roots by Lemma \ref{lhsrcomparison}. This says also that if $J$ is an irreducible subset of $I$, $J$ is not of type $\typeA$ and $I\in\calI(\be)$ then $J\in\calI(\be)$, i.e. $\calI(\be)$ is an ideal in the set of subsystems of type not $\typeA$. 
Further we introduce the set $\calF(\be)$ of fundamental weights $\om_{f(I)}$ for $I\in\calI(\be)$.

We will identify an element $I$ of $\calI(\be)$ with its type if this will not raise any confusion. For example if $\be=\theta$ for type $\typeC_4$, then $\calI(\be)=\{\typeC_4>\typeC_3>\typeB_2\}$. The only possible ambiguity is for $\Phi$ of type $\typeE$ that contains two different subsystems of type $\typeD_5$; in this case we denote by $\typeDL$ the subsystem corresponding to $\{1,2,3,4,5\}$ and by $\typeDR$ the one corresponding to $\{2,3,4,5,6\}$. The next lemma proves a property of $\calI(\be)$.
\begin{lemma}\label{lCDE}
Let $\Phi$ be an irreducible root system not of type $\typeA$ or $\tpG_2$
and let $\be=\sum_{i\in[\ell]}b_i\al_i$ with $\be\geq\zeta$. Suppose that $\hts(\be)\leq2$ and that $[\ell]\not\in\calI(\be)$. Then there exists $I\in\calI(\be)$ such that $\om_{f(I)}\in\supp^+(\be)$.
\end{lemma}
\begin{proof}
Before we treat the various types we want to make some general remarks. Notice first that $b_1,\ldots,b_\ell\geq1$ since $\be\geq\zeta$. We write also $\be$, as an element of the weight lattice, as
$\grb=\sum_{i\in[\ell]}a_i\om_i$.

Suppose now we have fixed the ideal $\calI=\calI(\be)$, this result in some simple conditions on the coefficient $b_1,\ldots,b_\ell$, also, $\calI$ determines the set $\calF=\calF(\be)$ (that sometimes we consider as a subset of $[\ell]$). Given such conditions on $b_1,\ldots,b_\ell$, the thesis is equivalent to show that the system
$$
\begin{array}{c}
\hts(\be)\leq 2\textrm{ and }a_j\leq 0\textrm{ for all }j\in\calF
\end{array}
$$
has no solution in $\grb$. 
Moreover notice that $\hts(\be)\leq 2$ is equivalent to $\sum_{j\in J}a_j\leq 2$ for all $J\subset[\ell]$ and to $\sum_{j\in J}a_j\leq 2$ for all $J\subset[\ell]\setminus\calF$.

Now we want deduce from the system above a new set of inequalities that will be useful in some case.
First notice that the system $a_j\leq 0$ for all $j\in\calF$ is clearly equivalent to the condition 
$\sum_{j\in\calF}x_j a_j\leq 0$ for all $(x_j)_{j\in\calF}\in(\R^+)^{|\calF|}$. Let 
$C=(c_{i,j})_{i,j\in[\ell]}$ be the Cartan matrix of $\De$, where 
$c_{i,j}=\langle \al_i,\al_j^\vee\rangle$, also let $C_\calF$ be the submatrix associated to $\calF$ 
and let $D_\calF\doteq C_\calF^{-1}=(d_{i,j})_{i,j\in\calF}$. The entries of $D_\calF$ are non negative, 
hence setting $x_j=\sum_{j\in\calF}d_{i,j}y_j$ we have $x_j\geq 0$  for all $i\in\calF$ provided that 
$y_j\geq 0$ for all $j\in\calF$. So we can write
$$
\sum_{j\in\calF}x_j a_j=\sum_{i\in\calF}b_i y_i+\sum_{i\not\in\calF}\big(\sum_{h,j\in\calF}c_{i,j}d_{j,h}y_h\big)b_i\leq0
$$
for every $(y_i)_{i\in\calF}\in(\R^+)^{|\calF|}$. In particular, if we let $y_i=\de_{i,t}$ we find
\begin{equation}\label{eq:systeminv}
b_t+\sum_{i\not\in\calF}\big(\sum_{j\in\calF}c_{i,j}d_{j,t}\big)b_i\leq 0\tag{$**$}
\end{equation}
for all $t\in\calF$. We notice also that
this system of inequalities is not equivalent 
to the original system $a_j\leq 0$ for all $j\in\calF$. 

In the following we fix $\grb$ and we set $\calI\doteq\calI(\grb)$ and $\calF\doteq\calF(\grb)$
and we conclude the proof with a case by case analysis.

\emph{Type} $\typeB_\ell$. In this case we observe that $\zeta = \theta = \gra_1 + \dots + \gra_\ell$,
hence there is nothing to prove.

\emph{Type} $\typeC_\ell$. For this type
$$
\calI=\{\tpC_h>\tpC_{h-1}>\cdots>\tpC_3>\tpB_2\}
$$
for some $h\geq3$ or $\calI=\{\tpB_2\}$ and we set $h=2$. The hypothesis $[\ell]\not\in\calI$ imposes $h<\ell$. Set $t\doteq\ell-h+1$. We have $\calF=\{\om_{t+1},\om_{t+2},\ldots,\om_{\ell-1},\om_\ell\}$. Using $\tpC_h\in\calI$, $\tpC_{h+1}\not\in\calI$ we have $b_t=1$, $b_{t+1}\geq2$. Notice that 
$a_{t+1}+a_{t+2}+\cdots+a_\ell=-b_t+b_{t+1}=b_{t+1}-1\geq1$, hence one of the integers $a_{t+1},a_{t+2},\ldots,a_\ell$, say $a_i$, must be positive. So $\om_i\in\calF\cap\supp^+(\be)$.

\emph{Type} $\typeD_\ell$. For this type $\calI=\varnothing$ or
$$
\calI=\{\tpD_h>\tpD_{h-1}>\cdots>\tpD_4\}
$$
for some $h\geq4$. Suppose $\calI\neq\varnothing$. As in the preceding case the hypothesis $[\ell]\not\in\calI$ gives $h<\ell$.  Set $t\doteq\ell-h+1$. We have $\calF=\{\om_{t+1},\om_{t+2},\ldots,\om_{\ell-2}\}$. Using $\tpD_h\in\calI$, $\tpD_{h+1}\not\in\calI$ we have $b_t=1$, $b_{t+1}\geq2$. Notice that 
$a_{\ell-1}+a_\ell =2(b_\ell+b_{\ell-1}-b_{\ell-2})$. Suppose that this is not positive. Then $a_{t+1}+a_{t+2}+\cdots+a_{\ell-2}=-b_t+b_{t+1}+b_{\ell-2}-b_{\ell-1}-b_\ell\geq1+b_{\ell-2}-b_{\ell-1}-b_\ell\geq1$, hence one of the integers $a_{t+1},a_{t+2},\ldots,a_{\ell-2}$, say $a_i$, must be positive. So $\om_i\in\calF\cap\supp^+(\be)$. On the other hand if $a_{\ell-1}+a_\ell$ is positive we have $a_{\ell-1}+a_\ell\geq2$. Now consider $a_1+a_2+\cdots+a_{t-1}=b_1+b_{t-1}-b_t=b_1+b_{t-1}-1\geq1$. Hence $\hts(\be)\geq a_1+a_2+\cdots+a_{t-1}+a_{\ell-1}+a_\ell\geq3$, and this is impossible since we assume $\hts(\be)\leq2$.

Now suppose $\calI=\varnothing$. Then $b_{\ell-2}=1$, hence
$a_{\ell-1}+a_\ell=2(b_\ell+b_{\ell-1}-b_{\ell-2})\geq2$. We have also $a_1+a_2+\cdots+a_{\ell-3}=b_1+b_{\ell-3}-b_{\ell-2}\geq1$. So $\hts(\be)\geq a_1+a_2+\cdots+a_{\ell-3}+a_{\ell-1}+a_\ell\geq3$ and this is impossible.

\emph{Type} $\typeE_6$. For this type $\calI$ is an ideal of the following poset
$$
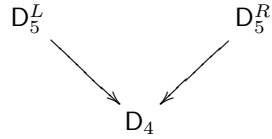
\xymatrix{
\tpDL\ar[dr] & & \tpDR\ar[dl]\\
 & \tpD_4\\
}
$$
So we have four possibilities, up to symmetries.\\
\emph{Case 1}: $\calI=\varnothing$. We have $\tpD_4\not\in\calI$, hence $b_4=1$. So $\hts(\be)\geq a_1+a_2+a_3+a_5+a_6=b_1+2b_2+b_3-3b_4+b_5+b_6=b_1+2b_2+b_3+b_5+b_6-3\geq3$; contrary to $\hts(\be)\leq 2.$\\
\emph{Case 2}: $\calI=\{\tpD_4\}$. We have $\tpDL,\tpDR\not\in\calI$, hence $b_3=b_5=1$, whereas $D_4\in\calI$ implies $b_4\geq 2$. So $\hts(\be)\geq a_1+a_2+a_4+a_6=2b_1+b_2+b_4+2b_6-2b_3-2b_5=2b_1+b_2+b_4+2b_6-4\geq 3$. So also this case is impossible.\\
\emph{Case 3}: $\calI=\{\tpDL,\tpD_4\}$, hence $\calF=\{\om_3,\om_4\}$. We have $\tpDR\not\in\calI$, so $b_5=1$. Also $2\geq\hts(\be)\geq a_1+a_2+a_3+a_4+a_6=b_1+b_2+2b_6-2\geq2$, hence $b_1=b_2=b_6=1$. We find $a_3+a_4=-b_1-b_2+b_3+b_4-1=b_3+b_4-3>0$. So $a_3>0$ or $a_4>0$ and this finishes the proof.\\
\emph{Case 4}: $\calI=\{\tpDL,\tpDR,\tpD_4\}$, hence $\calF=\{\om_3,\om_4,\om_5\}$. We have $b_3,b_4,b_5\geq2$, since $\tpDL,\tpDR\in\calI$, and either $b_4=2$ or $b_2=1$, since $\tpE_6\not\in\calI$. Suppose $b_4=2$. We have $2\geq\hts(\be)=a_1+a_2+a_3+a_5+a_6=b_1+2b_2+b_3+b_5-6\geq 2$, hence $b_1=b_2=1$, $b_3=b_5=2$. We find $a_3=2b_3-b_1-b_4=1$ (and $a_5=1$ too).

Now suppose $b_2=1$ and $b_4\geq 3$. If $a_3$ or $a_4$ or $a_5$ is positive we have finished, so suppose $a_3=2b_3-b_4-b_1\leq0$, $a_4=2b_4-b_3-b_5-1\leq0$ and $a_5=2b_5-b_4-b_6\leq0$, that is $2b_4\leq b_3+b_5+1$, $2b_3\leq b_1+b_4$ and $2b_5\leq b_4+b_6$. We have also $2\geq\hts(\be)\geq a_1+a_3+a_4+a_5+a_6=b_1+b_6-1$, hence $b_1+b_6\leq3$. We get $4b_4\leq 2b_3+2b_5+2\leq b_1+2b_4+b_6+2$, so $2b_4\leq b_1+b_6+2\leq 5$. This gives $b_4=1$ that is impossible.

\emph{Type} $\typeE_7$. For this type $\calI$ is an ideal of the following poset
$$
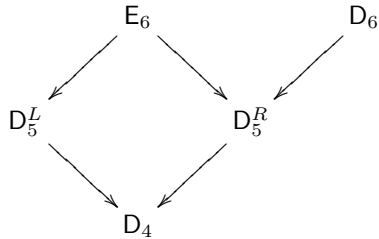
\xymatrix{
 & \tpE_6\ar[dl]\ar[dr] & & \tpD_6\ar[dl]\\
\tpDL\ar[dr] & & \tpDR\ar[dl]\\
 & \tpD_4\\
}
$$
The proof for the ideals that appear also in type $\tpE_6$ is systematically simpler than the one seen for that type; so we check only the cases of the new four ideals.\\
\emph{Case 1}: $\calI=\{\tpD_6,\tpDR,\tpD_4\}$, hence $\calF=\{\om_4,\om_5,\om_6\}$. Using $\tpDL\not\in\calI$ and $\tpD_6\in\calI$ we find $b_3=1$ and $b_4,b_5,b_6\geq 2$. Consider the inequalities
$$
\begin{array}{l}
2 \geq \hts(\be) \geq a_1+a_2+a_4+a_5+a_6+a_7\\
2 \geq \hts(\be) \geq a_1+a_4+a_5+a_6\\
2 \geq \hts(\be) \geq a_1+a_5\\
0 \geq a_6\\
\end{array}
$$
that we can write as
$$
\begin{array}{l}
4 \geq 2b_1+b_2+b_7\\
4+b_2+b_7 \geq 2b_1+b_4+b_6\\
3+b_4+b_6 \geq 2b_1+2b_5\\
b_5+b_7 \geq 2b_6.\\
\end{array}
$$
The first inequality gives $b_1=b_2=b_7=1$; so substituting these in the second one gives $b_4=b_6=2$ that in the third one gives $b_5=2$. Hence the last inequality becomes $3\geq 4$. So this case is impossible.\\
\emph{Case 2}: $\calI=\{\tpD_6,\tpDR,\tpDL,\tpD_4\}$, 
hence $\calF=\{\om_3,\om_4,\om_5,\om_6\}$. Since $\tpDL,\tpD_6\in\calI$ we get $b_3,b_4,b_5,b_6\geq 2$, whereas $\tpE_6\not\in\calI$ gives $b_4=2$ or $b_2=1$. Suppose first $b_4=2$ and consider the inequalities
$$
\begin{array}{l}
0\geq a_3\\
2\geq \hts(\be) \geq a_1+a_2+a_5+a_6+a_7\\
\end{array}
$$
that we can write as
$$
\begin{array}{l}
b_1 \geq 2b_3-2\\
6 \geq 2b_1-b_3+2b_2+b_5+b_7.\\
\end{array}
$$
Substituting the first in the second inequality we find $6\geq 3b_3+2b_2+b_5+b_7-4\geq 7$. This is impossible.\\
So we can suppose $b_4\geq 3$ and $b_2=1$. Consider the inequalities $2\geq\hts(\be)\geq a_1+a_7$ and the inequalities (\ref{eq:systeminv}) for $b_3,b_4,b_6$; we get the system
$$
\begin{array}{l}
2b_1+2b_7 \leq 2+b_3+b_6\\
b_3 \leq \frac{4}{5}b_1+\frac{1}{5}b_7+\frac{3}{5}\\
b_4 \leq \frac{3}{5}b_1+\frac{2}{5}b_7+\frac{6}{5}\\
b_6 \leq \frac{1}{5}b_1+\frac{4}{5}b_7+\frac{2}{5}.
\end{array}
$$
Adding up the second and the forth inequalities and substituting in the first one we get $2b_1+2b_7\leq 3+b_1+b_7$ i.e. $b_1+b_3\leq 3$. So, using the third inequality, we find $b_4\leq\frac{3}{5}(b_1+b_7)+\frac{6}{5}-\frac{1}{5}b_7\leq 3-\frac{1}{5}b_7<3$, contrary to the assumption $b_4\geq 3$.\\
\emph{Case 3}: $\calI=\{\tpE_6,\tpDL,\tpDR,\tpD_4\}$, hence $\calF=\{\om_2,\om_3,\om_4,\om_5\}$. From $\tpD_6\not\in\calI$ and $\tpE_6\in\calI$ we find $b_6=1$, $b_2,b_3,b_5\geq 2$ and $b_4\geq 3$. Consider the inequality $2\geq\hts(\be)\geq a_1+a_7$ and the inequality (\ref{eq:systeminv}) for $b_3$, we get the system
$$
\begin{array}{l}
2b_1+2b_7\leq 3+b_3,\\
b_3\leq b_1+\frac{1}{2},\textrm{ namely }b_3\leq b_1.
\end{array}
$$
Substituting the second inequality in the first one we get $b_3+2b_7\leq 3$, that is impossible since $b_3\geq 2$ and $b_7\geq 1$.\\
\emph{Case 4}: $\calI=\{\tpD_4,\tpDL,\tpDR,\tpE_6,\tpD_6\}$, hence $\calF=\{\om_2,\om_3,\om_4,\om_5,\om_6\}$. Using $\tpD_6,\tpE_6\in\calI$ we find $b_2,b_3,b_5,b_6\geq 2$ and $b_4\geq 3$, whereas using $\tpE_7\not\in\calI$ we find $b_1=1$ or $b_3=2$ or $b_4=3$ or $b_5=2$. First notice that the inequality $a_2\leq 0$ gives $b_4\geq 2b_2\geq 4$, hence $b_4=3$ is not possible.\\
If $b_3=2$ consider the inequality $2\geq\hts(\be)\geq a_1+a_2+a_4+a_5+a_6+a_7$ and the inequality (\ref{eq:systeminv}) for $b_4$, we get
$$
\begin{array}{l}
6\geq 2b_1+b_2+b_7\\
b_4\leq\frac{3}{2}b_1+b_7.
\end{array}
$$
From the first inequality we find $b_1=1$ and $b_7\leq 2$, so, substituting in the second inequality we find $b_4<4$. Hence the system has no solution since, as proved above, $b_4\geq 4$.\\
If $b_3\geq 3$ and $b_1=1$ consider the inequality $2\geq\hts(\be)\geq a_2+a_3+a_4+a_5+a_6+a_7$ and the inequality (\ref{eq:systeminv}) for $b_4$, we get
$$
\begin{array}{l}
b_2+b_3+b_7\leq 3+b_4\\
b_4\leq\frac{3}{2}+b_7,\textrm{ namely }b_4\leq 1+b_7.
\end{array}
$$
Substituting the second inequality in the first one we find $5\leq b_2+b_3\leq 4$ that is impossible.\\
If $b_3\geq 3$, $b_1\geq 2$ and $b_5=2$, consider the inequality $a_2\leq 0$ and $2\geq\hts(\be)\geq a_1+a_3+a_4+a_6+a_7=b_1+b_2+(b_4-2b_2)+b_6+b_7-4\geq 3$, hence the system above has no solution.

\emph{Type} $\typeE_8$. For this type $\calI$ is an ideal of the following poset
$$
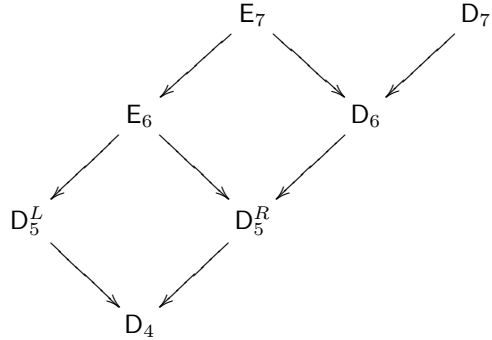
\xymatrix{
 & & \tpE_7\ar[dl]\ar[dr] & & \tpD_7\ar[dl]\\
 & \tpE_6\ar[dl]\ar[dr] & & \tpD_6\ar[dl]\\
\tpDL\ar[dr] & & \tpDR\ar[dl]\\
 & \tpD_4\\
}
$$
The proof for the ideals that already appear for $\tpE_7$ is simpler that the one for that type; so we check the remaining five ideals.\\
\emph{Case 1}: $\calI=\{\tpD_7,\tpD_6,\tpDR,\tpD_4\}$, 
hence $\calF=\{\om_4,\om_5,\om_6,\om_7\}$. Since $\tpDL\not\in\calI$ we have 
$b_3=1$, whereas $\tpD_7\in\calI$ gives $b_4,b_5,b_6,b_7\geq 2$. Using 
$2\geq\hts(\be)\geq a_1+a_2+a_4+a_5+a_6+a_7+a_8=2b_1+b_2+b_8-2$ we find 
$b_1=b_2=b_8=1$. So $a_4+a_5+a_6+a_7=b_4+b_7-b_2-b_3-b_8=b_4+b_7-3\geq 1$, 
hence at least one of the $a_4,a_5,a_6,a_7$ is positive.\\
\emph{Case 2}: $\calI=\{\tpD_7,\tpD_6,\tpDR,\tpDL,\tpD_4\}$, hence $\calF=\{\om_3,\om_4,\om_5,\om_6,\om_7\}$. Using $\tpD_7,\tpDL\in\calI$ we find $b_3,b_4,b_5,b_6,b_7\geq 2$, whereas $\tpE_6\not\in\calI$ gives $b_4=2$ or $b_2=1$. Suppose $b_4=2$. The inequality $2\geq\hts(\be)\geq a_1+a_2+a_3+a_5+a_6+a_7+a_8=b_1+b_3+2b_2+b_5+b_8-6\geq 2$ imposes $b_1=b_2=b_8=1$, $b_3=b_5=2$. So $a_3=2b_3-b_1-b_4=1$.\\
Now suppose $b_4\geq 3$ and $b_2=1$. Using $2\geq\hts(\be)\geq a_1+a_3+a_4+a_5+a_6+a_7+a_8=b_1+b_8-1$ we find $b_1+b_8\leq 3$. The inequality (\ref{eq:systeminv}) for $b_3$ reads
$$
b_3\leq\frac{5}{6}b_1+\frac{1}{6}b_8+\frac{2}{3}=\frac{5}{6}(b_1+b_8)+\frac{2}{3}-\frac{2}{3}b_8\leq\frac{19}{6}-\frac{2}{3}b_8\leq\frac{15}{6}<3
$$
hence $b_3=2$. In the same way the inequality (\ref{eq:systeminv}) for $b_5$ gives
$$
b_5\leq\frac{1}{2}(b_1+b_8)+1\leq\frac{5}{2}<3
$$
hence $b_5=2$. Thus $a_4=2b_4-b_2-b_3-b_5\geq 1$.\\
\emph{Case 3}: $\calI=\{\tpD_7,\tpD_6,\tpE_6,\tpDL,\tpDR,\tpD_4\}$, hence 
$\calF=\{\om_2,\om_3,\om_4,\om_5,\om_6,\om_7\}$. We have $b_2,b_3,b_5,b_6,b_7\geq 2$ and 
$b_4\geq 3$; also $\tpE_7\not\in\calI$ implies $b_1=1$ or $b_3=2$ or $b_4=3$ or $b_5=2$. 
Notice that $a_2\leq 0$ imposes $b_4\geq 2b_2\geq 4$, so $b_4=3$ is not possible. 
First suppose $b_3=2$. Using $2\geq\hts(\be)\geq a_1+a_2+a_4+a_5+a_6+a_7+a_8=2b_1+b_2+b_8-4$ we find 
$6\geq 2b_1+b_2+b_8$, hence $b_1=1$, $b_2\leq 3$ and $b_8\leq 2$. The inequality (\ref{eq:systeminv}) for 
$b_5$ gives
$$
b_5\leq\frac{3}{2}b_1+b_8\leq\frac{7}{2}
$$
so $b_5\leq 3$. Moreover $a_4=2b_4-b_3-b_2-b_5\leq 0$ gives $b_2=3$ and $b_4=4$, hence $a_2=2b_2-b_4=2>0$.\\
Now suppose $b_3\geq 3$ and $b_1=1$. Consider the inequality $2\geq\hts(\be)\geq a_2+a_3+a_4+a_5+a_6+a_7+a_8=b_2+b_3+b_8-b_4-1$, that is $b_2+b_3+b_8\leq b_4+3$, and combine it with the inequality (\ref{eq:systeminv}) for $b_4$, that is $b_4\leq 2+b_8$; we get $b_2+b_3\leq 5$. Hence $b_2=2$, $b_3=3$. Now we have $2\geq\hts(\be)\geq a_3+a_4+a_5+a_6+a_7+a_8=b_3+b_8-b_1-b_2=b_8$, hence $b_4\leq 2+b_8\leq 4$. So $b_4=4$ and we finish noting that $a_3=2b_3-b_1-b_4=1$.\\
Now suppose $b_3\geq3$, $b_1\geq 2$ and $b_5=2$. The inequality $2\geq\hts(\be)\geq a_1+a_2+a_3+a_4+a_6+a_7+a_8=b_1+b_2+b_6+b_8-4\geq 3$ shows that this case is not possible.\\
\emph{Case 4}: $\calI=\{\tpE_7,\tpD_6,\tpE_6,\tpDL,\tpDR,\tpD_4\}$, hence $\calF=\{\om_1,\om_2,\om_3,\om_4,\om_5,\om_6\}$. The condition $\tpE_7\in\calI$ gives $b_1\geq 2$, $b_2\geq 2$, $b_3\geq 3$, $b_4\geq 4$, $b_5\geq 3$ and $b_6\geq 2$, whereas $\tpD_7\not\in\calI$ gives $b_7=1$. This case is ruled out by inequality (\ref{eq:systeminv}) for $b_1$ (for example) giving $b_1\leq\frac{2}{3}$.\\
\emph{Case 5}: $\calI=\{\tpD_7,\tpE_7,\tpD_6,\tpE_6,\tpDL,\tpDR,\tpD_4\}$, hence 
$\calF=\{\om_1,\om_2,\om_3,\om_4,\om_5,\om_6,\om_7\}$. 
Using the condition $a_i\leq 0$ for all $i\in\calF$, we can write $\be=a\om_8-\eta$ 
for some integer $a$ and some dominant weight $\eta$ with support contained in $\calF$. 
Observe also 
that $\be\geq\zeta$ implies $\grb \notin -\grL^+$, hence $2\geq\hts(\be)=a$ imposes $a=1$ or $a=2$. 
Now notice that 
$\zeta=\sum_{i=1}^8\al_i=\om_1+\om_2-\om_4+\om_8$; also, the weight lattice coincides with the root 
lattice so we can write $\eta=\sum_{i=1}^8 c_i\al_i$ for some non negative integers $c_1,\ldots,c_8$. 
Suppose $a=1$. Writing the inequality $\be\geq\zeta$ in terms of the simple roots we get
$$
\eta=\sum_{i=1}^8 c_i\al_i\leq-\om_1-\om_2+\om_4=\al_1+2\al_2+3\al_3+5\al_4+4\al_5+3\al_6+2\al_7+\al_8.
$$
So $c_1\leq 1$, hence $\eta=0$ using the expression of the fundamental weights in terms of the simple roots. We find $\be=\om_8$ that implies $\tpE_8\in\calI$ contrary to our hypothesis.\\
If we suppose $a=2$ we have
$$
\eta=\sum_{i=1}^8 c_i\al_i\leq-\om_1-\om_2+\om_4+\om_8=3\al_1+5\al_2+7\al_3+11\al_4+9\al_5+7\al_6+5\al_7+3\al_8.
$$
But the unique fundamental weight with coefficient of $\al_1$ less or equal to $3$ is $\om_8$, hence $\eta=0$ since $\supp\eta\subset\calF$. We conclude as above that $\tpE_8\in\calI$ contrary to the hypothesis.

\emph{Type } $\typeF_4$: We observe first that 
$\grb \geq \zeta \geq \theta_{\typeB_3} = \gra_1+\gra_2+\gra_3$. Hence we have only two possibilities:
$\calI=\{\tpB_2,\tpB_3\}$ or $\calI=\{\tpB_2,\tpB_3,\tpC_3\}$.
In the first case $\calF=\{\om_1,\om_2\}$ implies $a_1,a_2 \leq 0$ and $\tpC_3 \notin \calI$ implies
$b_3=1$; so $a_3=2 - b_4 - 2b_2 \leq 0$.
In the second case $\calF=\{\om_1,\om_2,\om_3\}$ hence $a_1,a_2, a_3 \leq 0$.
In both case we have $\grb=a\om_4 - \eta$ with $\eta$ a dominant weight supported on $\{1,2,3\}$, also notice
that $2\geq \hts (\grb) = a$.
Now using $\grb \geq \zeta$ we find
$$
\eta \leq (a-1)\gra_1 + (2a-1)\gra_2 + (3a-1)\gra_3 + (2a-1) \gra_4.
$$
Looking at the coefficient of $\gra_1$ and arguing as in the last case for type $\tpE_8$, we 
deduce $\eta =0$. But this is impossible since it would give $\grb \geq \theta=\om_4$.
\end{proof}

A remark on the conclusion of Lemma \ref{lCDE}. 
Suppose that we can apply the lemma on $\be=\la+\mu-\nu$ for some 
$\la$, $\mu$ and $\nu$ dominant weights. We find an irreducible subset 
$I$ of $[\ell]$ such that $\be\geq\theta_I$ and $\om_{f(I)}\in\supp^+(\be)$. 
So $\om_{f(I)}\in\supp(\la)\cup\supp(\mu)$, hence $\la-\theta_I\in\La^+$ or 
$\mu-\theta_I\in\La^+$. So $(\la,\mu,\nu)$ is not a \low triple.

We are now ready to prove the main result of this section. First, we introduce a bit of notation. If $\la$ covers $\mu$ and $\la-\mu=\theta_I$, we write $\xymatrix@1{\la\ar[r]^I & \mu}$ (see Proposition \ref{pcovering}). A \emph{covering diagram} for a weight $\la\in\La^+$ is the direct graph whose vertexes are the dominant weights $\eta\leq\la$ and whose arrows are the covering relations. We see an example. Let $\Phi$ be of type $\typeB_3$ and let $\la=2\om_3$, then its covering diagram is
$$
\xymatrix@1{
2\om_3\ar[r]^{\{\ell\}} & \om_2\ar[r]^{\typeB_2} & \om_1\ar[r]^{\typeB_3} & 0
}
$$
since $2\om_3-\om_2=\al_\ell$, $\om_2-\om_1=\theta_{\{2,3\}}$ and $\om_1=\theta$. Observe that we use the identification of subsets of $[\ell]$ with their types as explained above Lemma \ref{lCDE}. Also notice that it is easy to write down the covering diagram of a given dominant weight: use Proposition \ref{pcovering} and the remark after the table of highest short roots. Further notice that, given a covering diagram
$$
\xymatrix@1{
\la_1\ar[r]^{I_1} & \la_2\ar[r]^{I_2} & \cdots\ar[r]^{I_r} & \la_{r+1}
}
$$
we have $\la_1-\la_{r+1}\geq\zeta$ if and only if $\cup_{i=1}^r I_i=[\ell]$.
\begin{proof}[Proof of Theorem~B]
We show first that any triple of the form $(\la,-w_0\la,0)$ with $\la$ minuscule is a \low triple. 
$\la$ being minuscule also $-w_0\la$ is minuscule and the property i) in Definition \ref{dlowtriple} is obvious. We have also $w_0\la\leq\la$ since $w_0\la$ is in the orbit $W\cdot\la$. Hence $\la-w_0\la=\be$ where $\be$ is some nonnegative linear combination of simple roots, say $\be=\sum_{\al\in\De}n_\al\al$. We want to show that $\supp_\Phi(\be)=\De$ proving property ii) in Definition \ref{dlowtriple}. Suppose that this is not the case and let $\ga$ be a simple root with $\ga\not\in\supp_\Phi(\De)$ adjacent in the Dynkin diagram to a simple root $\ga'\in\supp_\Phi(\be)$. Then in $\sum_{\al\in\De}n_\al\al$ the fundamental weight $\om_\ga$ dual to $\ga^\vee$ appears with a negative coefficient; this is not possible.

Now we see the other implication. Suppose $(\la,\mu,\nu)$ is a low triple and let
$\be=\la+\mu-\nu$. We want to show that $\la$, $\mu$ are minuscule, dual to each other and that $\nu=0$. 

We prove first that $\grl$ and $\mu$ are minuscule if $\Phi$ is not of type $\tpG_2$. 
Observe that, by Lemma \ref{lheight}, $ht^+_0(\grb) \leq 2$.
Observe also that $\grb \geq \theta$: indeed, if $\Phi$ is of type $\typeA$, this is clear by
$\grb \geq \zeta = \theta$ and, if $\Phi$ is not of type $\typeA$, it 
follows by the remark after the proof of Lemma \ref{lCDE}. 
Assume now that $\grl$ is not minuscule. Then there exists $I$ such that $\grl-\theta_I \in \grL^+$ and also 
$\grb \geq \theta \geq \theta_I$, so by Lemma \ref{luseful} $(\grl,\mu,\nu)$ is not a low a triple.
Similarly we prove that $\mu$ must be minuscule.

We perform now a case by case analysis assuming $\grl$ and $\mu$ minuscule.

\emph{Type} $\typeA_\ell$. 
In this case $\la$ and $\mu$ are fundamental weights, 
say $\la=\om_i$, $\mu=\om_j$ for some $i,j\in[\ell]$. 
Since $\la+\mu-\nu\geq\zeta$ the covering diagram of $\la+\mu$ must be of the following form
$$
\xymatrix@1{
\om_i+\om_j\ar[r]^(.4){\typeA(i,j)} & \om_{i-1}+\om_{j+1}\ar[r]^(.6){\typeA(i-1,j+1)} & \cdots \ar[r]^{\typeA(2,\ell-1)} & \om_1+\om_\ell\ar[r]^(.6){\typeA(1,\ell)} & 0
}
$$
where $\tpA(h,k)=\{h,h+1,\ldots,k\}$ and $\nu=0$. So $i+j=\ell+1$ proving our claim.

\emph{Type} $\typeB_\ell$. 
In this case we have $\la=\mu=\om_\ell$ that is the unique minuscule weights. 
We have the following covering diagram
$$
\xymatrix@1{
2\om_\ell\ar[r]^{\{\ell\}} & \om_{\ell-1}\ar[r]^{\tpB_2} & \om_{\ell-2}\ar[r]^{\tpB_3} & \cdots\ar[r]^{\tpB_{\ell-1}} & \om_1\ar[r]^{\tpB_\ell} & 0.
}
$$
Since $\la+\mu-\nu\geq\zeta$ we have $\nu=0$.

\emph{Type} $\typeC_\ell$. In this case we have that $\la=\mu=\om_1$. Thus the covering diagram
$$
\xymatrix@1{
2\om_1\ar[r]^{\{1\}} & \om_2\ar[r]^{\tpC_\ell} & 0.
}
$$
Since $\la+\mu-\nu\geq\zeta$ we have $\nu=0$. This proves our claim.

\emph{Type} $\typeD_\ell$. 
The possibilities, up to symmetry of the Dynkin diagram, are: 1) $\la=\mu=\om_1$, 2) $\la=\om_1$, $\mu=\om_\ell$, 3) $\la=\om_{\ell-1}$, $\mu=\om_\ell$ and 4) $\la=\mu=\om_\ell$. In each case we use $\la+\mu-\nu\geq\zeta$ to show $\nu=0$ and $\la$ and $\mu$ dual to each other.\\
\emph{Case 1}: we have
$$
\xymatrix@1{
2\om_1\ar[r]^{\{1\}} & \om_2\ar[r]^{\tpD_\ell} & 0.
}
$$
So the unique \low triple for this case is $(\om_1,\om_1,0)$.\\
\emph{Case 2}: we have
$$
\xymatrix@1{
\om_1+\om_\ell\ar[r]^(.6){\tpA^D_\ell} & \om_{\ell-1}
}
$$
where $\tpA^D_\ell=\{1,2,\ldots,\ell-2,\ell\}$. This shows that there is no \low triple in this case.\\
\emph{Case 3}: suppose $\ell$ odd, we have
$$
\xymatrix@1{
\om_{\ell-1}+\om_\ell\ar[r]^(.6){\tpA_3^R} & \om_{\ell-3}\ar[r]^{\tpD_5} & \om_{\ell-5}\ar[r]^{\tpD_7} & \cdots\ar[r]^{\tpD_{\ell-2}} & \om_2\ar[r]^{\tpD_\ell} & 0
}
$$
where $\tpA_3^R=\{\ell-1,\ell-2,\ell\}$. So the unique \low triple is $(\om_{\ell-1},\om_\ell,0)$. Now suppose $\ell$ even, we have
$$
\xymatrix@1{
\om_{\ell-1}+\om_\ell\ar[r]^(.6){\tpA_3^R} & \om_{\ell-3}\ar[r]^{\tpD_5} & \om_{\ell-5}\ar[r]^{\tpD_7} & \cdots\ar[r]^{\tpD_{\ell-3}} & \om_3\ar[r]^{\tpD_{\ell-1}} & \om_1.
}
$$
Hence there is no \low triple.\\
\emph{Case 4}: suppose $\ell$ odd, we have
$$
\xymatrix@1{
2\om_\ell\ar[r]^(.44){\{\ell\}} & \om_{\ell-2}\ar[r]^{\tpD_4} & \om_{\ell-4}\ar[r]^{\tpD_6} & \cdots\ar[r]^{\tpD_{\ell-3}} & \om_3\ar[r]^{\tpD_{\ell-1}} & \om_1.
}
$$
Hence there is no \low triple. Now suppose $\ell$ even, we have
$$
\xymatrix@1{
2\om_\ell\ar[r]^(.44){\{\ell\}} & \om_{\ell-2}\ar[r]^{\tpD_4} & \om_{\ell-4}\ar[r]^{\tpD_6} & \cdots\ar[r]^{\tpD_{\ell-2}} & \om_2\ar[r]^{\tpD_\ell} & 0.
}
$$
So $(\om_\ell,\om_\ell,0)$ is the unique \low triple.

\emph{Type} $\tpE$. As above we have that $\la$ and $\mu$ are minuscule weights. So we have the following cases up to symmetries: 1) $\Phi$ of type $\tpE_6$ and $\la=\mu=\om_1$, 2) $\Phi$ of type $\tpE_6$ and $\la=\om_1$, $\mu=\om_6$ and 3) $\Phi$ of type $\tpE_7$ and $\la=\mu=\om_7$. In each case we use $\la+\mu-\nu\geq\zeta$ to show that $\nu=0$ and $\la$ and $\mu$ are dual to each other.\\
\emph{Case } 1: we have
$$
\xymatrix@1{
2\om_1\ar[r]^{\{1\}} & \om_3\ar[r]^\tpDL & \om_6
}
$$
and so there is no \low triple in this case.\\
\emph{Case } 2: we have
$$
\xymatrix@1{
\om_1+\om_6\ar[r]^(.6){\tpA_5} & \om_2\ar[r]^{\tpE_6} & 0
}
$$
where $\tpA_5=\{1,3,4,5,6\}$. This says that the unique \low triple is $(\om_1,\om_6,0)$.\\
\emph{Case } 3: we have
$$
\xymatrix@1{
2\om_7\ar[r]^{\{7\}} & \om_6\ar[r]^{\tpD_6} & \om_1\ar[r]^{\tpE_7} & 0. 
}
$$
So the unique \low triple is $(\om_7,\om_7,0)$.

\emph{Type} $\tpF_4$. In this case there are not minuscule weights, hence there is nothing to prove.

\emph{Type} $\tpG_2$. We apply Lemma \ref{luseful} with $\grg=\zeta= \gra_1+\gra_2 =-\om_1+\om_2$ and 
$\grg=\gra_1$ to obtain $\grl=\mu=\om_1$.
So we have the following covering diagram
$$
\xymatrix@1{
2\om_1\ar[r]^{\{1\}} & \om_2\ar[r]^{\zeta} & \om_1\ar[r]^{\tpG_2} & 0. 
}
$$
Since $\nu\leq \grl + \mu -\zeta$ we find $\nu= \om_1$ or $\nu = 0$.
In both cases if $\grl'=0$ we have $\nu \leq \grl' + \mu$ and $\grl'\leq \grl$ 
against the first condition for a low triple.
\end{proof}

\subsection{Non reduced root system}\mbox{}\\
We still need to treat the case of the irreducible non reduced root system of type $\typeBC_\ell$. We begin with a word of explanation about this type. We think of $\Phi$ as the union of $\typeB_\ell$, with square root lengths $1,2$, and $\typeC_\ell$ with square root lengths $2,4$. The base $\De=\{\al_1,\ldots,\al_\ell\}$ for $\Phi$ is the same of $\typeB_\ell$ with indexing from \cite{B}; also $\al_\ell$ is the unique simple root such that $2\al_\ell\in\Phi$. One can define the set $\La$ of integral weights $\la$ requiring that $\langle \la,\coal\rangle\in\Z$ for all $\al\in\Phi$. As one can easily see, a weight $\la$ is integral if and only if $\langle\la,\coal_i\rangle\in\Z$ for $i=1,\ldots,\ell-1$ and $\langle \la,(2\al_\ell)^\vee\rangle\in\Z$. So the fundamental weights $\om_1,\ldots,\om_\ell$ are those of $\typeC_\ell$. Further we order $\La$ as usual defining $\mu\leq\la$ if and only if $\la-\mu\in R^+$. So the same definition of \low triples seen for reduced root system applies also to this situation.

We begin noticing that for this type there is no minuscule weight, i.e. the weight lattice coincides with the root lattice.
In the following proposition we see that Theorem~B is still valid also for non reduced root systems.
\begin{proposition}\label{pnonreducedlowtriple}
If $\Phi$ is an irreducible non reduced root system then there is no \low triple.
\end{proposition}
\begin{proof}
Suppose that $(\la,\mu,\nu)$ is a \low triple. Lemma \ref{luseful} holds also for type $\typeBC$. Then by the expression $\sum_{i=h}^\ell\al_i=-\om_{h-1}+\om_h$, we have that $h\not\in\supp\la\cup\supp\mu$, finishing the proof.
\end{proof}

\subsection{Enlarged weight lattice}\mbox{}\\
In order to treat the case of exceptional complete symmetric variety we need a slight generalization of Theorem~B in case of $\Phi$ of type $\typeBC_\ell$.

Let us realize such type in an euclidean space $E$ and consider $E\simeq E\times0\hookrightarrow E\times\R$. Let $\om_1,\ldots,\om_\ell\in E\times0$ be the fundamental weights of $\Phi$ and take two linear independent vectors $u,v\in E\times\R\setminus E\times0$ such that $u+v=\om_\ell$. Consider the following cones $\La^+=\langle\om_1,\ldots,\om_\ell\rangle_\N$, $P^+=\langle\om_1,\ldots,\om_{\ell-1},u,v\rangle_\N$ and the related lattices $\La=\La^+_\Z$ and $P=P^+_\Z$. We can order $P$ declaring $\mu\leq\la$ if and only if $\la-\mu\in R^+$, where $R^+\subset E\times0$ is the cone generated by the simple roots. We define low triples for weights in $P^+$ using the same definition as in \ref{dlowtriple}.
\begin{lemma}\label{lenlargedbelow}
Let $\nu\in P^+$ and $\la=\la'+au\in P^+$ for some $\la'\in\La^+$ and $a\in\N$. Then $\nu\leq\la$ if and only if $\nu=\nu'+au$ for some $\La^+\ni\nu'\leq\la'$. The same holds for $\la=\la'+av$.
\end{lemma}
\begin{proof}
This is clear from the definition of the order since $R\subset\La$.
\end{proof}
We want to show that
\begin{proposition}\label{pexclowtriple}
Let $\la,\mu,\nu$ be three weights in $P^+$. If $(\la,\mu,\nu)$ is a low triple then, up to symmetry, $\la=au$, $\mu=bv$ with $a,b\geq1$. If $a\geq b$ then $\nu=\nu'+(a-b)\om_\al$ with $\La^+\ni\nu'\leq b\om_\ell$. Moreover if $\la=u$, $\mu=v$ then $\nu=0$.
\end{proposition}
\begin{proof}
Notice that $\sum_{i=h}^\ell\al_i=-\om_{h-1}+\om_h$ for $h=1,\ldots,\ell$, as already seen for type $\typeBC_\ell$. This shows that $\la=a_1u+a_2v$, $\mu=b_1u+b_2v$. Also $a_1a_2=0$, since otherwise $\la=\om_\ell+\la'$ for some $\la'\in P^+$ and this is not possible. This shows that $\la,\mu\in\N u\cup\N v$. 
Further if $\la=au$, $\mu=bu$ then, using Lemma \ref{lenlargedbelow} above, 
$\la+\mu=(a+b)u$ is a minimal element in $P^+$; this is not possible since $P^+\ni\nu<\la+\mu$. So we must have $\la=au$, $\mu=bv$ up to symmetry.

Now suppose that $a\geq b$. Then $\la+\mu=au+bv=b\om_\ell+(b-a)u$ and the stated form of $\nu$ follows by Lemma \ref{lenlargedbelow} above. Finally notice that $(u,v,\nu)$ is a low triple in $P^+$ if and only if $(2\om_\ell,2\om_\ell,\nu)$ is a low triple for type $\typeB_\ell$, as follows directly since we are using the simple roots of $\typeB_\ell$ and $2\om_\ell$ is the fundamental weight dual to $\coal_\ell$. So $\nu=0$ by Theorem~B.
\end{proof}


\vskip 0.5cm
\scriptsize{
Rocco Chiriv\`\i,\\
Universit\`a di Pisa,\\
Dipartimento di Matematica ``Leonida Tonelli''\\
via Buonarroti n.~2, 56127 Pisa, Italy,\\
e-mail: {\tt chirivi@dm.unipi.it}\\
\indent and\\
Andrea Maffei,\\
Universit\`a di Roma ``La Sapienza",\\
Dipartimento di Matematica ``Guido Castelnuovo",\\
Piazzale Aldo Moro n.~5, 00185 Roma, Italy,\\
e-mail: {\tt amaffei@mat.uniroma1.it}
}
\end{document}